\documentclass[12pt,a4paper]{article}

\usepackage[british]{babel}

\usepackage[a4paper,top=2cm,bottom=2cm,left=2.5cm,right=2.5cm,marginparwidth=1.75cm]{geometry}

\usepackage[style=apa, backend=biber]{biblatex} % APA 7th edition style citations using biblatex
\DeclareLanguageMapping{british}{british-apa} % Set language mapping
\DeclareFieldFormat[article]{volume}{\apanum{#1}} % Format volume number

\usepackage{amsmath}
\usepackage{graphicx}
\usepackage[colorlinks=true, allcolors=blue]{hyperref}
\usepackage{hyperref}
\usepackage[title]{appendix}
\usepackage{multirow}
\usepackage{mathrsfs}
\usepackage{amsfonts}
\usepackage{nicematrix}
\usepackage{amssymb}
\usepackage{amsmath}
\usepackage{booktabs} % For \toprule, \midrule, \botrule
\usepackage{caption}  % For \caption
\usepackage{threeparttable} % For table footnotes
\usepackage{algorithm}
\usepackage{algorithmicx}
\usepackage{algpseudocode}
\usepackage{listings}
\usepackage{enumitem}
\usepackage{chngcntr}
\usepackage{booktabs}
\usepackage{lipsum}
\usepackage{subcaption}
\usepackage{authblk}
\usepackage[T1]{fontenc}    % Font encoding
\usepackage{csquotes}       % Include csquotes
\usepackage{diagbox}
\usepackage{upgreek}

\usepackage{setspace}
\usepackage{titlesec}
\titleformat{\section} % Redefine section numbering format
  {\normalfont\Large\bfseries}{\thesection.}{1em}{}
  
\usepackage{lineno} % Add the lineno package

\rightlinenumbers % Right-align line numbers

\linenumbers % Enable line numbering

\usepackage{float}   % for customizing caption position
\usepackage{caption} % for customizing caption format


\title{Asymptotics of higher-order conditional tail moments for convolution-equivalently distributed losses}
\author[1]{Zhangting Chen}
\author[2]{Bingjie Wang}
\author[1,*]{Dongya Cheng}
\affil[1]{\small School of Mathematical Sciences, Soochow University, Suzhou 215006, P.R. China}
\affil[2]{School of Mathematical Science and LPMC, Nankai University, Tianjin 300071, P.R. China}
\affil[*]{Corresponding author: \texttt{dycheng@suda.edu.cn}}

\date{}  % Remove date

\begin{document}
\newtheorem{lemma}{Lemma}
\newtheorem{pron}{Proposition}
\newtheorem{assu}{Assumption}
\newenvironment{assubis}[1]
  {\renewcommand{\theassu}{\ref{#1}$'$}%
   \addtocounter{assu}{-1}%
   \begin{assu}}
  {\end{assu}}
\newtheorem{re}{Remark}
\newtheorem{thm}{Theorem}
\newtheorem{co}{Corollary}
\newtheorem{exam}{Example}
\newtheorem{defin}{Definition}
\newtheorem{proof}{Proof}
\newcommand{\la}{\frac{1}{\lambda}}
\newcommand{\sectemul}{\arabic{section}}
\renewcommand{\theequation}{\sectemul.\arabic{equation}}
\renewcommand{\thepron}{\sectemul.\arabic{pron}}
\renewcommand{\thelemma}{\sectemul.\arabic{lemma}}
\renewcommand{\there}{\sectemul.\arabic{re}}
\renewcommand{\thethm}{\sectemul.\arabic{thm}}
\renewcommand{\theco}{\sectemul.\arabic{co}}
\renewcommand{\theexam}{\sectemul.\arabic{exam}}
\renewcommand{\thedefin}{\sectemul.\arabic{defin}}
\renewcommand{\theassu}{\sectemul.\arabic{assu}}
\maketitle

\begin{abstract}
This paper investigates the asymptotic behavior of higher-order conditional tail moments, which quantify the contribution of individual losses in the event of systemic collapse. The study is conducted within a framework comprising two investment portfolios experiencing dependent losses that follow convolution-equivalent distributions. The main results are encapsulated in two theorems: one addressing light-tailed losses with convolution-equivalent distributions and the other focusing on heavy-tailed losses with regularly varying distributions. Both results reveal that the asymptotic behavior remains robust regardless of the strength of dependence. Additionally, numerical simulations are performed under specific scenarios to validate the theoretical results. 
%The simulation results confirm the accuracy and reliability of the proposed asymptotic approximations.
%Abstracts must be able to stand alone and cannot contain citations to the paper’s references, equations, etc. An abstract must consist of a single paragraph and be concise. Because of online formatting, abstracts must appear as plain as possible. Three to six keywords must be included. Each keyword should not exceed three words. %\lipsum[1]
\end{abstract}

\textbf{Keywords}: Asymptotics, Dependence, Light tail, Regular variation, Higher-order conditional tail moment.

{\bf{Mathematics Subject Classification (2020)}} Primary 62P05 $\cdot$ Secondary 91G70

%-------------------------------------------
% Paper Body
%-------------------------------------------
\section{Introduction}\setcounter{thm}{0}\setcounter{equation}{0}\setcounter{lemma}{0}\setcounter{co}{0}\setcounter{re}{0}

This paper considers two portfolios, where one portfolio (here we call it the investment portfolio) consists of $n$ ($n\geq 2$) assets (these assets could be stocks, bonds, real estate, commodities, or alternative investments), and the other portfolio (here we call it the hedging portfolio) consists of the corresponding $n$ hedging assets serving to offset or reduce the risk of the investment portfolio. For example,
\begin{itemize}
	\item A portfolio containing a stock can be hedged by stock index futures or put options against potential losses.
	\item For portfolios with international exposure (e.g., foreign stocks or bonds), then it might face the risk of currency fluctuations. In such cases, it can be hedged against currency risks by using foreign exchange derivatives, such as foreign exchange forwards, foreign exchange options, or currency exchange-traded Funds.
	\item In periods of macroeconomic uncertainty, stock market volatility, or inflation, investors often face significant risks, and gold and precious metals are common hedges against inflation and stock market declines.  
\end{itemize} 

In this paper, we denote the net losses of assets of the investment portfolio by $X_1$, $X_2$, $\cdots$, $X_n$, and the net losses of assets of the hedging portfolio by $Y_1$, $Y_2$, $\cdots$, $Y_n$, then the aggregate loss of the investment portfolio can be expressed by $S_n=\sum_{i=1}^nX_i$ and the aggregate loss of the hedging portfolio is formulated by $T_n=\sum_{i=1}^nY_i$. 

Risk measures are essential in finance because they provide a structured way to quantify, assess, and manage the risks faced by financial institutions, investors, and regulators. Understanding and managing risk is crucial for ensuring financial stability, making informed investment decisions, and achieving long-term goals while mitigating potential losses. Value at risk (VaR) is a typical risk measure. However, the 2008 Financial Crisis, the 2011 European Debt Crisis, and the Swiss Franc Peg Crisis are prime examples where VaR failed to predict the extent of the losses, for it did not account for extreme events or tail risks (rare but high-impact scenarios). Therefore, many researchers have intensively investigated various variants of VaR. The reader is referred to {Adrian and Brunnermeier (2016)} for conditional VaR (CoVaR) and $\Delta$CoVaR, to {Acharya et al. (2017)} for the marginal expected shortfall (MES for short, and under the present notational convention, it is formulated as $E(X_k|S_n>x)$) and systemic expected shortfall (SES).

Recently, {Ji et al. (2021)} introduced a risk measure called joint expected shortfall (JES), which integrates both expected shortfall (ES, for short, and under the present notational convention, it is formulated as $E(X_k|X_k>x)$) and MES and can be quantified as the expected loss for the individual when both the individual and the system are performing poorly. Motivated by {Ji et al. (2021)}, this paper considers the following quantity to fully address the impact of market uncertainty on portfolios. For $0\leq \zeta\leq 1$, $\beta>0$, and $1\leq k\leq n$, we define
\begin{align}\label{TM}
	\rho_k(x,y;\beta, \zeta)=E\left(X_k^{\beta}\Big\vert X_k>\zeta x, S_n>x, T_n>y\right),
\end{align}
where $x$ and $y$ represent thresholds, $\beta$ serves as the order of moments, and $\zeta$ quantifies the conservatism in the conditioning event and provides additional flexibility when the given information indicates that both $X_k$ and $S_n$ are large but not necessarily at the same scale. One can deduce from (\ref{TM}) that if we take $\zeta=0$, then $\rho_k(x,y;\beta, 0)$ reduces to a variation of MES, and if we take $\zeta=1$, $\rho_k(x,y;\beta, 1)$ becomes a variation of ES when $X_k$, $1\leq k\leq n$, are non-negative. Our main results include the asymptotic behavior of ES and MES since they are crucial in predicting the extent of the losses. From the perspective of the mathematical formula, $\rho_k(x,y;\beta, \zeta)$ is the $\beta-$order conditional expectation of the $k$th asset's loss of the investment portfolio when its value exceeds a certain threshold, meanwhile, the aggregate losses $S_n$ and $T_n$ exceed their respective thresholds, meaning both two portfolios collapse. Noticing that we can replace $X_k$ by $Y_k$ and $\zeta x$ by $\zeta y$ to reflect the contribution of the $k$th asset's loss to the hedging portfolio, but we will not do a specific study in this paper because they are perfectly symmetrical.

As is well documented in the literature, due to some technical difficulties, such as the mutual interplay among the losses, the complexity of the losses' distributions, and the magnitude of the order $\beta$, it is challenging to derive the close form of $\rho_k(x,y;\beta, \zeta)$, and thus it is difficult to quickly determine the specific value under the difficulties mentioned earlier which contradicts the need for rapid risk assessments and responses when a crisis occurs. Many scholars have sought to address this contradiction by focusing on deriving asymptotic formulas. As a matter of fact, under the current prudent insurance regulatory frameworks, this transition is not without justification, especially for the insurance industry. Some regulatory frameworks, for example, Solvency II, require that solvency capital should cover one-year losses with a confidence level of at least 99.5\%. In practice, many individual insurance companies select an even more conservative confidence level, as high as 99.9\%. In this regard, the asymptotic formula for (\ref{TM}) is rather reasonable.

Among those works seeking asymptotic formulas for risk measures, many of them assumed that the distributions of risks have regular variation or dominated variation, which means that the risks have extremely heavy tails. Some recent works are listed as follows. {Das and Fasen-Hartmann (2018)} derived the asymptotic behavior of both the marginal mean excess and the MES as the threshold goes to infinity, although the losses are asymptotically tail-independent in the model. {Liu and Yang (2021)} considered the asymptotic behavior of the systemic risk when the portfolio losses in the model are allowed to be regularly-varying-tailed losses, which are equipped with a wide range of dependence structures. {Chen and Liu (2022)} derived the asymptotic formulas of the SES and MES with long-tailed and dominatedly-varying-tailed losses. {Li (2022)} considered a dynamic systemic risk measure based on a multi-dimensional renewal risk model. When the losses possess pairwise asymptotical independence or multivariate regular variation, two asymptotic results of the systemic risk were respectively obtained in {Li (2022)}. The asymptotic behavior of the tail moment (TM) and the tail central moment (TCM) was considered in {Li (2023)} when the losses are pairwise asymptotically independent and regularly-varying-tailed. 

It should be noted that when the distributions of losses are max-stable distributions, some asymptotic formulas of risk measures are obtained. {Asimit and Li (2018)} investigated the asymptotics of the ES when the distributions of losses are max-stable. {Ji et al. (2021)} derived asymptotic formulas for the JES when the distributions of losses belong to the maximum domain of attraction of the Fr\'echet distribution, the Weibull distribution, or the Gumbel distribution. The dependence structure between losses in {Ji et al. (2021)} was described by a survival copula satisfying some regular variation properties. {Li (2023)} also obtained asymptotic results of the TM and the TCM when the losses belong to the maximum domain of attraction of the Gumbel distribution. 

Although the heavy-tailed distributions are an appealing choice to model the losses for their theoretical value and practical relevance in risk theory, the light-tailed distributions are also useful in insurance contexts. One can refer to {Gschl\"o\ss l and Czado (2007)} and {Wang and Hob\ae k Haff (2019)} for some related discussions. In the insurance context, modeling insurance losses by heavy-tailed distributions may cause a significant overestimation of solvency capital requirements, further leading to a higher cost of capital. For applications and research of light-tailed claims in risk theory, we refer the reader to {Cheng et al. (2002)}, {Tang and Tsitsiashvili (2003)}, {Yang and Yuen (2016)}, {Liu and Yang (2017)}, {Chen et al. (2023)}, among many others.

Recently, following the idea of {Li (2023)}, {Wang and Li (2024)} studied the asymptotic behavior of the TM when the losses follow convolution-equivalent or Gamma-like distributions. Although the distributions of losses are assumed to be exponential-like-tailed, which belong to the maximum domain of attraction of the Gumbel distribution, with the auxiliary function being a constant. Inspired by {Wang and Li (2024)}, we intend to derive the asymptotic formula of $\rho_k(x,y;\beta, \zeta)$ with losses' distributions belonging to the convolution-equivalent class, a broad family of distributions encompassing both heavy-tailed and light-tailed cases. 

Our model explicitly accounts for the dependence between hedging asset losses and investment asset losses. To this end, we assume that $\left\{(X_i, Y_i), i\geq 1\right\}$ is an independent and identically distributed (i.i.d.) random loss pair sequence, where the joint distribution of $(X_1, Y_1)$ is given in terms of the Farlie-Gumbel-Morgenstern (FGM) copula with parameter $\theta\in [-1,1]$, the reader is referred to {Nelsen (2006)} for a comprehensive discussion about copulas. Our work fills in some gaps in the asymptotic study of $\rho_k(x,y;\beta, \zeta)$ for dependent convolution-equivalently distributed losses. The lemmas in this paper have the potential to be applied in many other fields related to risk theory. Our main result is in the concise form of $C\times x^{\beta}$ for some $C>0$, which demonstrates that the strength of dependence does not influence the asymptotic behavior and assists us in calculating $\rho_k(x,y;\beta,\zeta)$ efficiently.
\section{Preliminaries and main results}\setcounter{thm}{0}\setcounter{equation}{0}\setcounter{lemma}{0}\setcounter{co}{0}\setcounter{re}{0}
Throughout the paper, unless otherwise noted, all limit relations hold as $x\to\infty$ or $x\wedge y\to\infty$. For two positive functions $f$ and $g$, we write $f\lesssim g$ if $\limsup f/g\leq 1$; write $f=O(g)$ if $\limsup f/g<\infty$; write $f=o(g)$ if $\limsup f/g=0$; write $f\sim g$ if $\lim f/g=1$. Let $f*g$ denote the convolution of $f$ and $g$, and $f^{*n}$ denote the $n$-fold of convolution with $f$ itself. For two real numbers $a$ and $b$, we write $a\vee b=\max\{a,b\}$, $a\wedge b=\min\{a,b\}$ and $a^{+}=a{\bf 1}(a\geq 0)$, $a^{-}=-a{\bf 1}(a< 0)$, where ${\bf 1}(A)$ represents the indicator function of set $A$. Moreover, if there is no possibility of confusion, we always write $\widehat{V}(\gamma)=\int_{-\infty}^{\infty}e^{\gamma x}V(\mathrm{d}x)$ for any distribution $V$ and any constant $\gamma\geq 0$.
\subsection{Distribution classes}
In the following, we shall introduce some distribution classes that will be used in the rest of this paper. For any distribution $V$ supported on $(-\infty,\infty)$, we denote its tail by $\overline{V}(x)=1-V(x)$ for all $x\in(-\infty,\infty)$. It is said to be regularly-varying-tailed with index $\alpha \geq 0$ and denoted by $V\in \mathscr{R}_{-\alpha}$ if its right tail is regularly varying, i.e.,
\begin{align*}
	\lim\frac{\overline{V}(xy)}{\overline{V}(x)}=y^{-\alpha}
\end{align*}
holds for all $y> 0$. The class of rapidly-varying-tailed distributions, $\mathscr{R}_{-\infty}$, is characterized by relation 
\begin{align*}
	\lim\frac{\overline{V}(xy)}{\overline{V}(x)}=0
\end{align*}
for any $y>1$. Notably, the class $\mathscr{R}_{-\infty}$ contains both heavy-tailed and light-tailed distributions such as the lognormal, the Weibull, and the exponential distributions. One closely related distribution class is the exponential-like-tailed distribution class, denoted by $\mathscr{L}(\gamma)$, for $\gamma\geq 0$. The distribution $V$ is said to belong to $\mathscr{L}(\gamma)$ for some $\gamma\geq 0$ if for every $y\in (-\infty,\infty)$, it holds that
\begin{align}\label{LG}
	\lim\frac{\overline{V}(x-y)}{\overline{V}(x)}=e^{\gamma y}.
\end{align}
It is called the uniform convergence theorem that for each $t>0$, the limit relation (\ref{LG}) holds uniformly for $y\leq t$ if $\gamma>0$, and uniformly for $-t\leq y\leq t$ if $\gamma=0$, see, for example, Pakes (2004). In addition, it is obvious that $V\in\mathscr{L}(\gamma)$ if and only if $\overline{V}\circ \log$ is regularly varying with index $\gamma$. By Proposition 2.2.3 of Bingham et al. (1987), if $V\in\mathscr{L}(\gamma)$ with $\gamma>0$, then for $0<\varepsilon<\gamma$ and $y>0$, there exists $x_0>0$ such that
\begin{align}\label{impor}
	\frac{\overline{V}\left(x+\frac{y}{\gamma}\right)}{\overline{V}(x)}\leq (1+\varepsilon)\exp\left\{-\left(1-\frac{\varepsilon}{\gamma}\right)y\right\}
\end{align} 
holds for $x\geq x_0$.

For $\gamma\geq 0$, a famous subclass of $\mathscr{L}(\gamma)$ is the convolution-equivalent distribution class, denoted by $\mathscr{S}(\gamma)$. If $V\in \mathscr{L}(\gamma)$, and the limit relation
\begin{align}\label{SG}
	\lim \frac{\overline{V*V}(x)}{\overline{V}(x)}=2d<\infty
\end{align}
holds for some $d>0$, then it is said to belong to $\mathscr{S}(\gamma)$, denoted by $V\in\mathscr{S}(\gamma)$. It is well known that the constant $d$ in (\ref{SG}) is equal to $\widehat{V}(\gamma)$, see {Rogozin (2000)} or {Foss and Korshunov (2007)} for related discussions and for $\gamma\geq 0$, $V\in\mathscr{S}(\gamma)$ if and only if $V^{+}\in\mathscr{S}(\gamma)$, see Corollary 2.1 (i) in {Pakes (2004)}. Remarkably, the distribution classes $\mathscr{L}(0)$ and $\mathscr{S}(0)$ retrieve the long-tailed distribution class $\mathscr{L}$ and the subexponential distribution class $\mathscr{S}$, respectively. The reader is referred to monographs such as {Embrechts et al. (1997)} and {Foss et al. (2011)} for some discussions about heavy-tailed distributions. We close this subsection by giving two assertions, which are implied respectively by Corollary 3.3.32 and Page 50 in {Embrechts et al. (1997)}  that 
\begin{align*}
	&\bullet \mathscr{S}(\gamma)\subset \mathscr{L}(\gamma)\subset \mathscr{R}_{-\infty},\quad {\rm for}\quad {\rm any}\quad \gamma>0;\\
	&\bullet \mathscr{R}_{-\alpha}\subset \mathscr{S}\subset \mathscr{L},\quad  {\rm for}\quad {\rm any}\quad \alpha\geq 0.
\end{align*}
\subsection{The FGM distribution}
Many scholars have proposed numerous dependence structures to describe various dependent phenomena. The reader is referred to {Chen and Yuen (2009)}, {Yang and Wang (2013)}, and {Li (2018)}, among many others, for some appealing and commonly used ones. As an efficient and widely used dependence structure, the FGM distribution is adopted in many papers and plays an important role in many fields of risk theory. Let $(X, Y)$ be a random pair whose marginal distributions are $F$ and $G$, respectively. It is said to follow the FGM distribution if for all $x, y\in (-\infty,\infty)$ the joint distribution is of the form
\begin{align}\label{joint}
	P(X\leq x, Y\leq y)={F}(x){G}(y)(1+\theta \overline{F}(x)\overline{G}(y))
\end{align}
with parameter $\theta\in[-1,1]$ governing the strength of the dependence. The following decomposition will be frequently used in the subsequent proofs. We observe that the joint survival function $\Pi$, which is defined by $\Pi(x,y)=P(X>x,Y>y)$ for all $x,y\in (-\infty,\infty)$, satisfies
\begin{align*}
	\Pi=(1+\theta)FG-\theta F^2G-\theta FG^2+\theta F^2G^2-F-G+1,
\end{align*}
or, equivalently,
\begin{align*}
	\Pi=(1+\theta)\overline{F}\ \overline{G}-\theta \overline{F}^2\overline{G}-\theta \overline{F}\ \overline{G}^2+\theta\overline{F}^2\overline{G}^2.
\end{align*}
The above relation suggests that if a random pair $(X, Y)$ follows a FGM distribution with parameter $\theta\in[-1,1]$, then it holds for any $x, y\in (-\infty, \infty)$ that
\begin{align}\label{decomp}
	\Pi(x,y)&=(1+\theta)P(X^{*}> x, Y^{*}> y)-\theta P(X^{*}_{\wedge}> x, Y^{*}> y)\nonumber\\
	&\quad -\theta P(X^{*}> x, Y^{*}_{\wedge}> y)+\theta P(X^{*}_{\wedge}> x, Y^{*}_{\wedge}> y),
\end{align}
where mutually independent random variables (r.v.s) $X^{*}$, $Y^{*}$, $X^{*}_{\wedge}$, and $Y^{*}_{\wedge}$ are distributed by $F$, $G$, $2F-F^2$, and $2G-G^2$, respectively. One can refer to {Liu and Yang (2017)} or {Chen et al. (2023)} for some related discussions.
%\section{The main result and numerical simulations}

\subsection{Main results}
After introducing some preliminaries about distribution classes and the FGM distribution, we will give the main results of this paper in the following. In the following theorem, we consider the case of light-tailed losses with convolution-equivalent distributions.
\begin{thm}\label{thm1}
Consider $\rho_k(x,y;\beta, \zeta)$ defined by (\ref{TM}). Let $\left\{(X, Y),(X_i, Y_i), 1\leq i\leq n\right\}$ be a sequence of i.i.d. random pairs with marginal distributions $F, G\in\mathscr{S}(\gamma)$, $\gamma\geq 0$, respectively. Let the joint distribution of $(X, Y)$ satisfy (\ref{joint}) with parameter $\theta\in [-1,1]$, then 
\begin{itemize}
	\item if $0<\zeta<1$, it holds for all $1\leq k\leq n$ and $\beta\in\mathbb{R}_{+}$ that
		\begin{align}\label{thm1res}
		\rho_k(x,y;\beta, \zeta)\sim x^{\beta};
		\end{align}
	\item if $\zeta=0$, it holds for all $1\leq k\leq n$ and $\beta\in\mathbb{N}_{+}$ that
	 \begin{align}\label{thm1res-1}
		\rho_k(x,y;\beta, \zeta)\sim \frac{x^{\beta}}{n}.
		\end{align}
\end{itemize}
\end{thm}
The next theorem considers the case of heavy-tailed losses with regularly varying distributions, where $F$ and $G$ are regularly varying.
\begin{thm}\label{thm2}
Consider $\rho_k(x,y;\beta, \zeta)$ defined by (\ref{TM}). Let $\left\{(X, Y),(X_i, Y_i), 1\leq i\leq n\right\}$ be a sequence of i.i.d. random pairs with marginal distributions $F, G\in\mathscr{R}_{-\alpha}$, $\alpha> 0$, respectively. Let the joint distribution of $(X, Y)$ satisfy (\ref{joint}) with parameter $\theta\in [-1,1]$, then 
\begin{itemize}
	\item if $0<\zeta\leq 1$, it holds for all $1\leq k\leq n$ and $\beta\in(0,\alpha)$ that 
		\begin{align}\label{thm2res}
	\rho_k(x,y;\beta, \zeta)\sim \frac{\alpha}{\alpha-\beta} x^{\beta};
		\end{align}
	\item if $\zeta=0$, $\alpha>1$, and $F(-x)=O\left(\overline{F}(x)\right)$, it holds for all $1\leq k\leq n$ and $\beta\in\mathbb{N}_{+}\cap (0,\alpha)$ that
		\begin{align}\label{thm2res-1}
	\rho_k(x,y;\beta, \zeta)\sim \frac{\alpha}{\alpha-\beta}\frac{x^{\beta}}{n}.
		\end{align}
\end{itemize}
\end{thm}
\begin{re}
	We observe that the strength of dependence (i.e., $\theta$) does not influence the asymptotic behavior of $\rho_k(x,y;\beta, \zeta)$. Furthermore, although the left-hand sides of (\ref{thm1res})-(\ref{thm2res-1}) depend on both $x$ and $y$, the asymptotic results are solely determined by $x$.
\end{re}
\begin{re}
	From the subsequent proof, it can be observed that in (\ref{thm1res-1}) and (\ref{thm2res-1}), if $X_i$, $1\leq i\leq n$, are non-negative, the condition on $\beta$ can be relaxed to allow non-integer values. This relaxation is particularly relevant because, to facilitate a more comprehensive comparison of the asymptotic behavior across different cases, $\beta$ is set as a real number in the numerical simulations presented in the next section.
\end{re}
\section{Numerical simulations}\setcounter{thm}{0}\setcounter{equation}{0}\setcounter{lemma}{0}\setcounter{co}{0}\setcounter{re}{0}
This section is designed to verify the accuracy of the relations (\ref{thm1res})-(\ref{thm2res-1}) when $\beta$, $\theta$, and $G$ vary, respectively. The following numerical simulations will set $k=1$, $n=3$, and $x=y$ to facilitate the following operations without losing rationality. The specific simulation methods in this paper are as follows.
\begin{itemize}
	\item Generating $N=10^6$ samples of i.i.d. random pairs $\left\{\left(X_{1,j}, Y_{1,j}\right),\left(X_{2,j}, Y_{2,j}\right),\left(X_{3,j}, Y_{3,j}\right)\right\}$, $j=1$, $2$, $\cdots$, $N$, where the joint distribution of each pair follows FGM copula with parameter $\theta$.
	\item Let $S_{3,j}=\sum_{i=1}^3X_{i,j}$ and $T_{3,j}=\sum_{i=1}^3Y_{i,j}$, then the empirical estimate of $\rho_1(x,x;\beta, \zeta)$ is constructed by
		\begin{align*}
			\widehat{\rho_1}(x,x;\beta, \zeta)=\frac{\sum_{j=1}^NX^{\beta}_{1,j}{\bf 1}_{\left(X_{1,j}>\zeta x, S_{3,j}>x, T_{3,j}>x\right)}}{\sum_{j=1}^N{\bf 1}_{\left(X_{1,j}>\zeta x, S_{3,j}>x, T_{3,j}>x\right)}},\quad 0\leq\zeta\leq 1.
		\end{align*}
	\item Calculating the asymptotic values of $\rho_1(x,x;\beta, \zeta)$ by relations (\ref{thm1res})-(\ref{thm2res-1}) for each $x$.
\end{itemize}

Before proceeding to the next step of the numerical simulations, it is necessary to briefly summarize the procedure of our simulations. We validate the asymptotic results of Theorems \ref{thm1}–\ref{thm2} via a two-stage numerical experiment. In the setting of Theorem \ref{thm1}, we respectively focus on the scenarios where $(\gamma,\zeta)=(\frac1{2},0)$ and $(\gamma,\zeta)=(\frac1{2},\frac1{3})$ while varying $\beta$, $\theta$, and $G$. For conditions of Theorem \ref{thm2}, we replicate the protocol.

\subsection{The case of $\gamma>0$}
We validate the accuracy of (\ref{thm1res}). In this scenario, we assume that both $F$ and $G$ are the inverse Gaussian distributions with the density functions $f_F(x)=\sqrt{\frac{1}{2\pi x^3}}\exp\left\{-\frac{(x-1)^2}{2x}\right\}$, and $f_G(x)=\sqrt{\frac{\nu}{2\pi x^3}}\exp\left\{-\frac{\nu(x-\mu)^2}{2\mu^2x}\right\}$, respectively, where $x>0$, $\mu>0$, and $\nu>0$. According to {Embrechts (1983)}, it can be verified  that $F\in \mathscr{S}\left(\frac1{2}\right)$ and $G\in \mathscr{S}\left(\frac{\nu}{2\mu^2}\right)$. To satisfy the conditions of Theorem \ref{thm1}, we consistently assume $\nu=\mu^2$ in the remainder of this section.

\begin{figure}[htbp] %H为当前位置，!htb为忽略美学标准，htbp为浮动图形
\centering %图片居中
\includegraphics[width=0.8\textwidth]{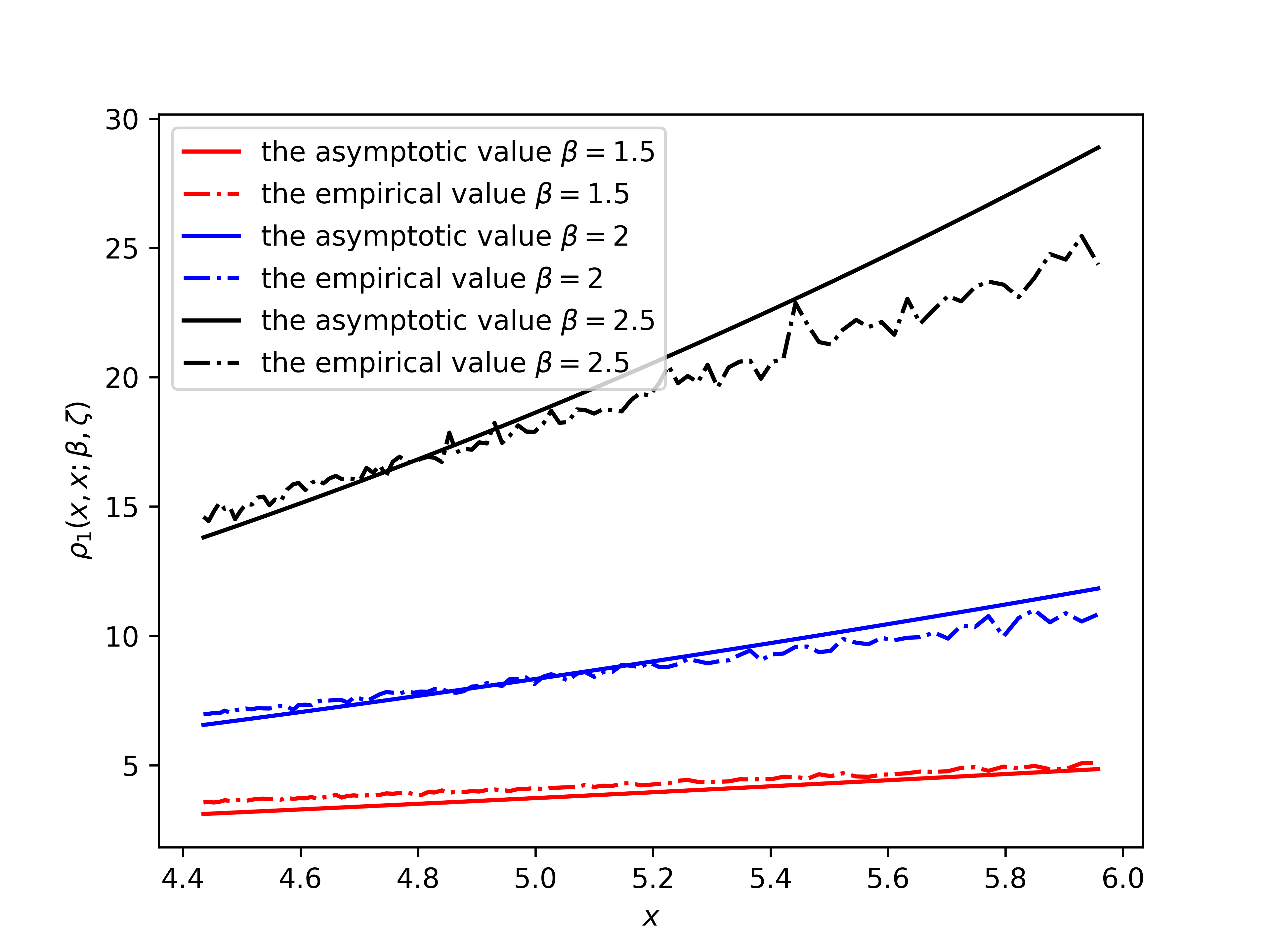} %插入图片，[]中设置图片大小，{}中是图片文件名
\captionsetup{font={footnotesize}}
\caption{Comparison between the asymptotic values for $\rho_1(x,x;\beta, \zeta)$ and its empirical values with $\theta=0.5$, $\zeta=0$, $\nu=1$, and $\mu=1$.} % 最终文档中希望显示的图片标题
\label{Fig.main4} %用于文内引用的标签
\end{figure}

\begin{figure}[htbp] %H为当前位置，!htb为忽略美学标准，htbp为浮动图形
\centering %图片居中
\includegraphics[width=0.8\textwidth]{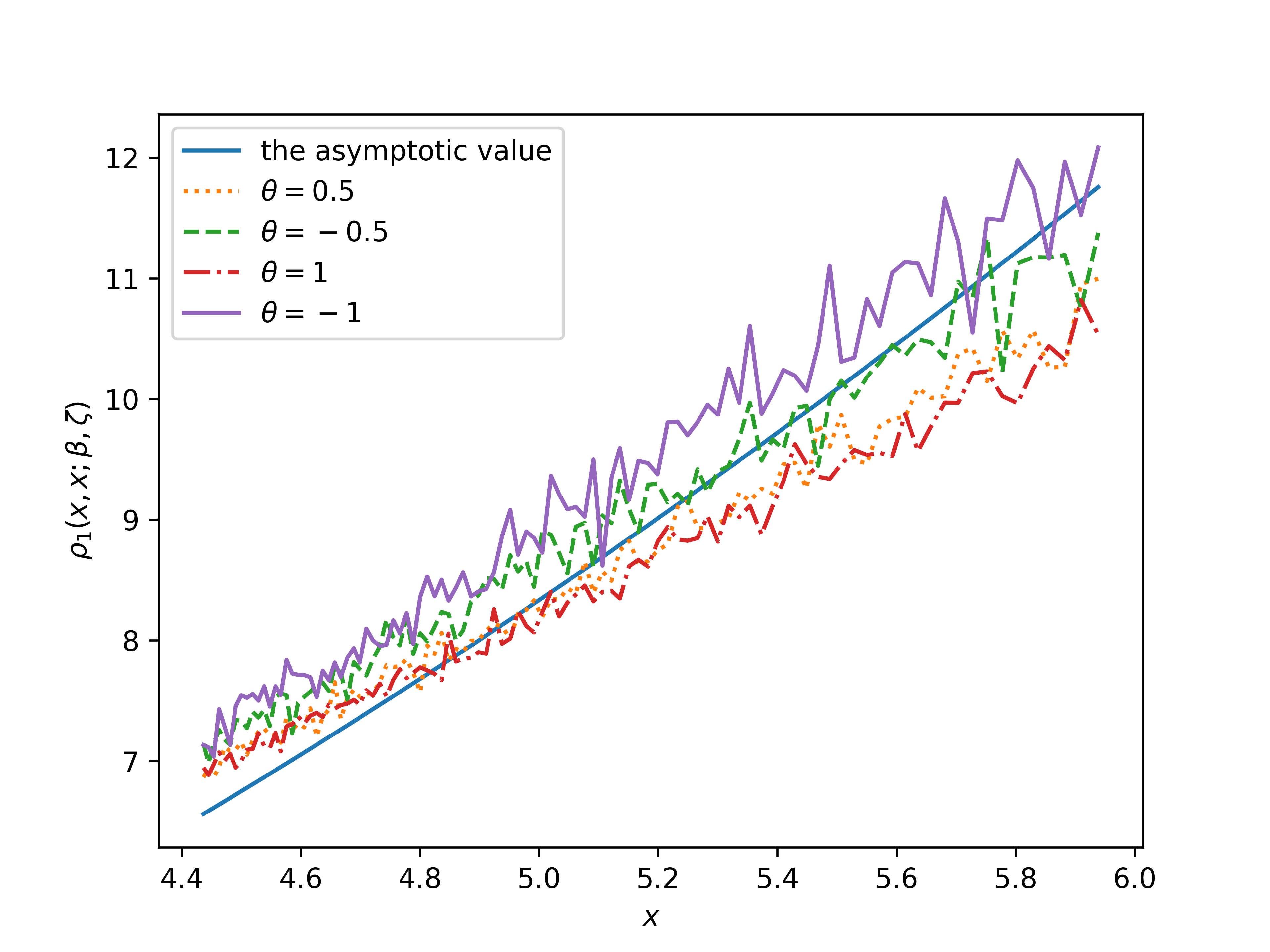} %插入图片，[]中设置图片大小，{}中是图片文件名
\captionsetup{font={footnotesize}}
\caption{Comparison between the asymptotic values for $\rho_1(x,x;\beta, \zeta)$ and its empirical values with $\beta=2$, $\zeta=0$, $\nu=1$, and $\mu=1$.} % 最终文档中希望显示的图片标题
\label{Fig.main5} %用于文内引用的标签
\end{figure}

\begin{figure}[htbp] %H为当前位置，!htb为忽略美学标准，htbp为浮动图形
\centering %图片居中
\includegraphics[width=0.8\textwidth]{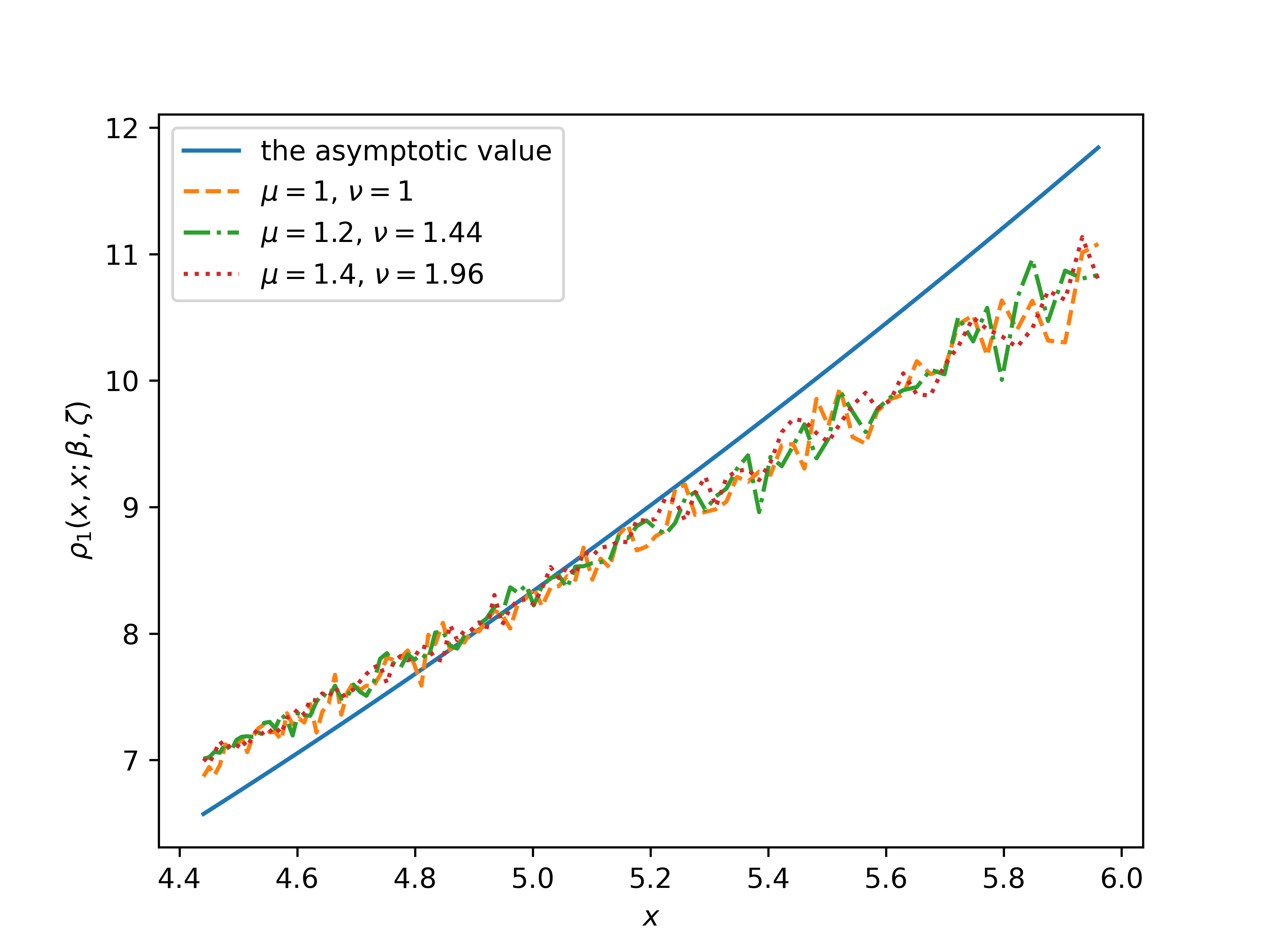} %插入图片，[]中设置图片大小，{}中是图片文件名
\captionsetup{font={footnotesize}}
\caption{Comparison between the asymptotic values for $\rho_1(x,x;\beta, \zeta)$ and its empirical values with $\beta=2$, $\theta=0.5$, $\zeta=0$, and $\gamma=\frac{1}{2}$.} % 最终文档中希望显示的图片标题
\label{Fig.main6} %用于文内引用的标签
\end{figure}

In Figures \ref{Fig.main4}-\ref{Fig.main6}, we set $\zeta=0$ to examine the impact of different parameters on the asymptotic estimates. In Figure \ref{Fig.main4}, we assume $\theta=0.5$, $\nu=1$, and $\mu=1$. The simulation results are presented in Figure \ref{Fig.main4}, indicating good performance of the asymptotic formula (\ref{thm1res}) but with fluctuations. Furthermore, as shown in Figure \ref{Fig.main4}, increasing $\beta$ leads to larger fluctuations in the error. This is to be expected, as higher powers tend to amplify the error.

In Figure \ref{Fig.main5}, we set $\nu=1$, $\mu=1$, and $\beta=2$ to investigate the accuracy of the asymptotic formula (\ref{thm1res}) with the variation of $\theta$. 
It can be deduced from Figure \ref{Fig.main5} that the dependence strength influences the empirical values, thereby affecting the accuracy of the estimates. Although the asymptotic formula (\ref{thm1res}) performs adequately when doing estimates, it is possible that a more precise asymptotic expression, involving the parameter $\theta$, could provide better estimates. Due to the inherent randomness in the simulation, we cannot make definitive conclusions at this stage, but this observation offers a valuable direction for future work.

Next, we investigate the accuracy of the asymptotic formula (\ref{thm1res}) by varying the distribution $G$ under the condition $\nu=\mu^2$. The corresponding results are plotted in Figure \ref{Fig.main6}. As shown in Figure \ref{Fig.main6}, the variation in distribution $G$ appears to have little influence on the empirical value of $\rho_1(x,x;\beta, \zeta)$. Our asymptotic result provides an accurate estimate of $\rho_1(x,x;\beta, \zeta)$ when $x$ falls within the range of 4.8 to 5.2. However, the asymptotic estimate tends to overestimate $\rho_1(x,x;\beta, \zeta)$ for $x>5.2$, while it underestimates the value for $x<4.8$. This behavior is expected, as the probability of a r.v. following an inverse Gaussian distribution exceeding 5.2 is less than 0.009, whereas the probability of it being below 4.8 exceeds 0.988. To achieve more accurate results, increasing the sample size $N$ would be necessary; however, this would come at the cost of significantly higher computational demands.

In Tables \ref{tab-1}-\ref{tab-3}, we will assume $\zeta=\frac1{3}$, while keeping all other settings consistent with the previous three studies. Detailed simulation results are provided in Tables \ref{tab-1}-\ref{tab-3} below for further reference.

\begin{table}[htbp]
\centering
\captionsetup{font={footnotesize}}
\caption{Ratios of asymptotic values to simulated values with $\theta=0.5$, $\zeta=\frac{1}{3}$, $\mu=1$, and $\nu=1$.}
\label{tab-1}
{\footnotesize
\resizebox{0.8\linewidth}{!}{
\begin{tabular}{cccccc}
\toprule[2pt]
\multicolumn{1}{c}{\multirow{2}{*}{$\beta$}} & \multicolumn{5}{c}{ratio} \\ \cline{2-6} 
\multicolumn{1}{c}{}& $x=2.0$ & $x=2.1$ & $x=2.2$ & $x=2.3$ & $x=2.4$ \\ \hline
{1}    & 0.9718  & 0.9436    & 1.0256 & 0.9878   & 1.0743 \\
{1.5}  & 0.9874  & 1.0108    & 0.9515 & 0.9955   & 0.9708 \\
{2}    & 1.0872  & 1.0393    & 0.9931 & 0.9591   & 0.9317 \\ 
\bottomrule[2pt]
\end{tabular}}}
\end{table}

\begin{table}[htbp]
\centering
\captionsetup{font={footnotesize}}
\caption{Ratios of asymptotic values to simulated values with $\beta=2$, $\zeta=\frac1{3}$, $\mu=1$, and $\nu=1$.}
\label{tab-2}
{\footnotesize
\resizebox{0.8\linewidth}{!}{
\begin{tabular}{cccccc}
\toprule[2pt]
\multicolumn{1}{c}{\multirow{2}{*}{$\theta$}} & \multicolumn{5}{c}{ratio} \\ \cline{2-6} 
\multicolumn{1}{c}{}& $x=2.0$ & $x=2.1$ & $x=2.2$ & $x=2.3$ & $x=2.4$ \\ \hline
{-1}   &1.0433   &0.9982     &0.9656 &0.9323   &0.9106 \\
%\textbf{-0.75}&1.1294   &0.9979     &0.9439 &0.9360   &0.9695 \\
{-0.5} &1.0591   &1.0133     &0.9758 &0.9425   &0.9181 \\
%\textbf{-0.25}&0.9124   &1.0309     &0.9435 &0.9834   &1.0803 \\
{0}    &1.0722   &1.0304     &0.9900 &0.9535   &0.9156 \\
%\textbf{0.25} &1.0338   &1.0962     &1.2063 &1.0118   &0.9835 \\
{0.5}  &1.0872   &1.0394     &0.9931 &0.9592   &0.9317 \\
%\textbf{0.75} &0.9489   &1.0267     &0.9454 &1.1633   &1.0705 \\
{1}    &1.1050   &1.0504     &1.0025 &0.9628   &0.9246 \\ 
\bottomrule[2pt]
\end{tabular}}}
\end{table}

\begin{table}[htbp]
\centering
\captionsetup{font={footnotesize}}
\caption{Ratios of asymptotic values to simulated values with $\beta=2$, $\theta=0.5$, and $\zeta=\frac1{3}$.}
\label{tab-3}
{\footnotesize
\resizebox{0.8\linewidth}{!}{
\begin{tabular}{cccccc}
\toprule[2pt]
\multicolumn{1}{c}{\multirow{2}{*}{($\mu$,$\nu$)}} & \multicolumn{5}{c}{ratio}\\ \cline{2-6} 
\multicolumn{1}{c}{}& $x=2.0$ & $x=2.1$ & $x=2.2$ & $x=2.3$ & $x=2.4$ \\ \hline
%\textbf{2.1}&1.0691   &1.2022     &0.8741 &1.1595   &1.0770 \\
%\textbf{2.3}&1.1294   &0.9979     &0.9439 &0.9360   &0.9695 \\
{(1,1)}       &1.0974   &1.0463     &1.0034 &0.9622   &0.9312 \\
{(1.2,1.44)}  &1.0642   &1.0218     &0.9785 &0.9492   &0.9248 \\
{(1.4,1.96)}  &1.0974   &1.0463     &1.0034 &0.9622   &0.9314 \\
%\textbf{3}  &1.0338   &1.0962     &1.2063 &1.0118   &0.9835 \\
\bottomrule[2pt]
\end{tabular}}}
\end{table}

% Please add the following required packages to your document preamble:
% \usepackage{multirow}
\begin{table}[htbp]
\centering
\captionsetup{font={footnotesize}}
\caption{Ratios of asymptotic values to simulated values with $\theta=0.5$, $\zeta=0$, and $\alpha=2.8$.}
\label{tab1}
{\footnotesize
\resizebox{0.8\linewidth}{!}{
\begin{tabular}{cccccc}
\toprule[2pt]
\multicolumn{1}{c}{\multirow{2}{*}{$\beta$}} & \multicolumn{5}{c}{ratio} \\ \cline{2-6} 
\multicolumn{1}{c}{}& $q=98.5\%$ & $q=98.7\%$ & $q=98.9\%$ & $q=99.1\%$ & $q=99.3\%$ \\ \hline
{1}    & 0.8820  & 0.8126    & 0.9837 & 0.9548   & 1.0147 \\
{1.5}  & 0.8306  & 1.0108    & 1.0239 & 0.8418   & 0.9957 \\
{2}    & 0.7436  & 1.2338    & 0.9531 & 1.0817   & 1.0321 \\ 
\bottomrule[2pt]
\end{tabular}}}
\end{table}

\begin{table}[htbp]
\centering
\captionsetup{font={footnotesize}}
\caption{Ratios of asymptotic values to simulated values with $\beta=2$, $\zeta=0$, and $\alpha=2.8$.}
\label{tab2}
{\footnotesize
\resizebox{0.8\linewidth}{!}{
\begin{tabular}{cccccc}
\toprule[2pt]
\multicolumn{1}{c}{\multirow{2}{*}{$\theta$}} & \multicolumn{5}{c}{ratio} \\ \cline{2-6} 
\multicolumn{1}{c}{}& $q=98.5\%$ & $q=98.7\%$ & $q=98.9\%$ & $q=99.1\%$ & $q=99.3\%$ \\ \hline
{-1}   &1.0691   &1.2022     &0.8741 &1.1595   &1.0770 \\
%\textbf{-0.75}&1.1294   &0.9979     &0.9439 &0.9360   &0.9695 \\
{-0.5} &0.9178   &1.1049     &0.9735 &1.2107   &1.0381 \\
%\textbf{-0.25}&0.9124   &1.0309     &0.9435 &0.9834   &1.0803 \\
{0}    &1.0580   &0.9763     &1.0692 &0.9612   &1.1104 \\
%\textbf{0.25} &1.0338   &1.0962     &1.2063 &1.0118   &0.9835 \\
{0.5}  &1.1683   &0.8082     &1.2044 &0.9407   &1.1039 \\
%\textbf{0.75} &0.9489   &1.0267     &0.9454 &1.1633   &1.0705 \\
{1}    &1.1356   &0.9523     &1.0449 &0.9496   &1.0043 \\ 
\bottomrule[2pt]
\end{tabular}}}
\end{table}

\begin{table}[htbp]
\centering
\captionsetup{font={footnotesize}}
\caption{Ratios of asymptotic values to simulated values with $\beta=2$, $\theta=0.5$, and $\zeta=0$.}
\label{tab3}
{\footnotesize
\resizebox{0.8\linewidth}{!}{
\begin{tabular}{cccccc}
\toprule[2pt]
\multicolumn{1}{c}{\multirow{2}{*}{$\alpha$}} & \multicolumn{5}{c}{ratio}\\ \cline{2-6} 
\multicolumn{1}{c}{}& $q=98.5\%$ & $q=98.7\%$ & $q=98.9\%$ & $q=99.1\%$ & $q=99.3\%$ \\ \hline
%\textbf{2.1}&1.0691   &1.2022     &0.8741 &1.1595   &1.0770 \\
%\textbf{2.3}&1.1294   &0.9979     &0.9439 &0.9360   &0.9695 \\
{2.5}&1.0647   &1.0630     &1.0194 &1.0943   &0.9735 \\
{3}  &1.1147   &1.0852     &1.0765 &1.0545   &1.0491 \\
{3.5}&1.1063   &1.0876     &1.0761 &0.9605   &1.0490 \\
%\textbf{3}  &1.0338   &1.0962     &1.2063 &1.0118   &0.9835 \\
\bottomrule[2pt]
\end{tabular}}}
\end{table}

Tables \ref{tab-1}-\ref{tab-3} imply that the majority of the ratios of the simulated and asymptotic values fall in the interval $[0.95, 1.05]$ when some involved parameters vary. 

\subsection{The case of $\gamma=0$}
In this scenario, set $F(x)=G(x)=1-\left(\frac{1}{x+1}\right)^{\alpha}$, $x\geq 0$, and $F,G\in\mathscr{R}_{-\alpha}\subset \mathscr{S}(0)$. To present the limiting relationship, let $x=y=F^{-1}(q)$, where $F^{-1}(q)=\inf\{x: F(x)\geq q\}$ for $q\in (0,1)$. It is obvious that $x=y=F^{-1}(q)\uparrow\infty$ as $q\uparrow 1$.  

\begin{figure}[htbp] %H为当前位置，!htb为忽略美学标准，htbp为浮动图形
\centering %图片居中
\includegraphics[width=0.8\textwidth]{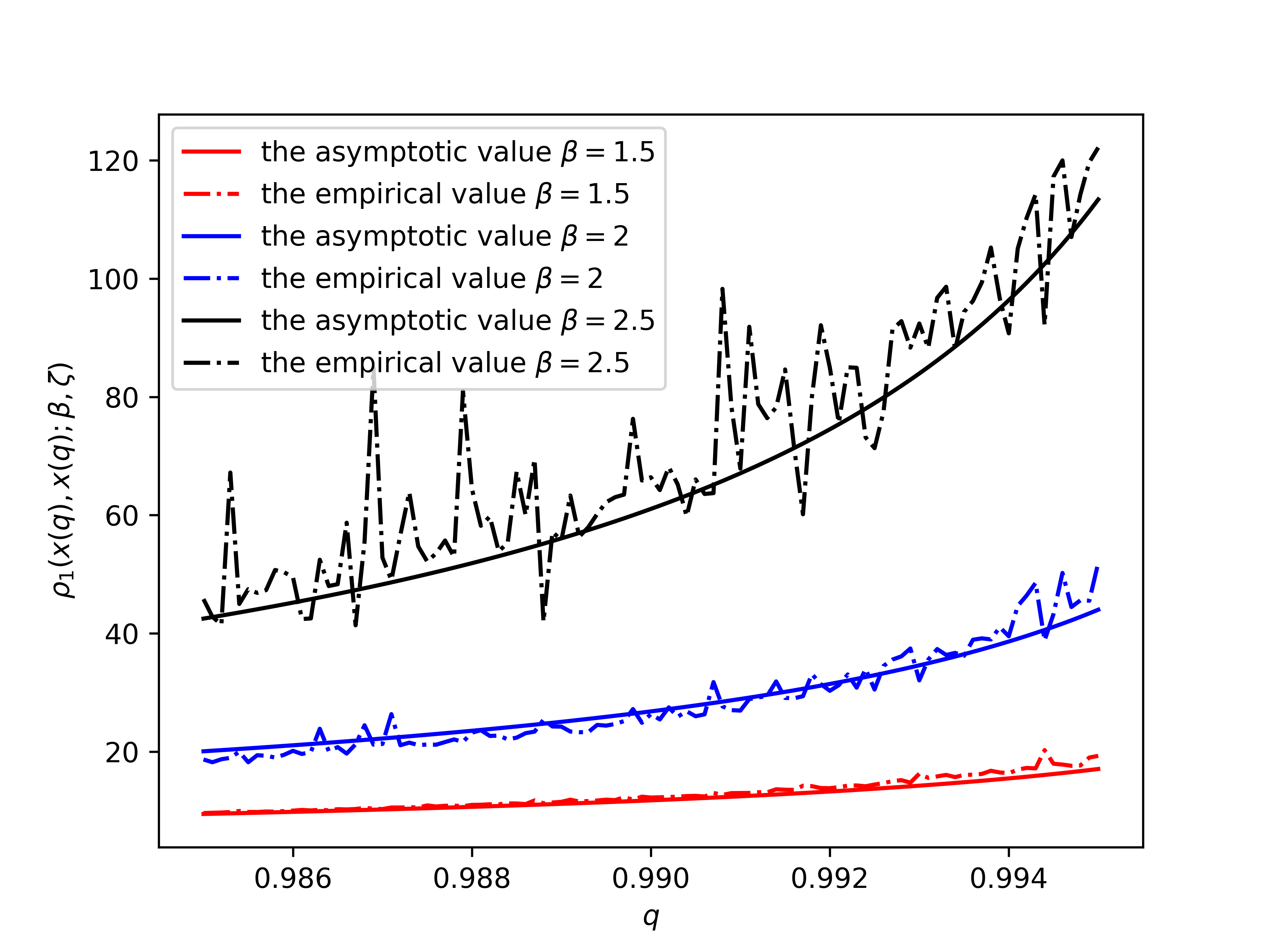} %插入图片，[]中设置图片大小，{}中是图片文件名
\captionsetup{font={footnotesize}}
\caption{Comparison between the asymptotic values for $\rho_1(x,x;\beta, \zeta)$  and its empirical values with $\theta=0.5$, $\zeta=\frac1{3}$, and $\alpha=2.8$.} % 最终文档中希望显示的图片标题
\label{Fig.main7} %用于文内引用的标签heavy_tailed_thtea
\end{figure}

\begin{figure}[htbp] %H为当前位置，!htb为忽略美学标准，htbp为浮动图形
\centering %图片居中
\includegraphics[width=0.8\textwidth]{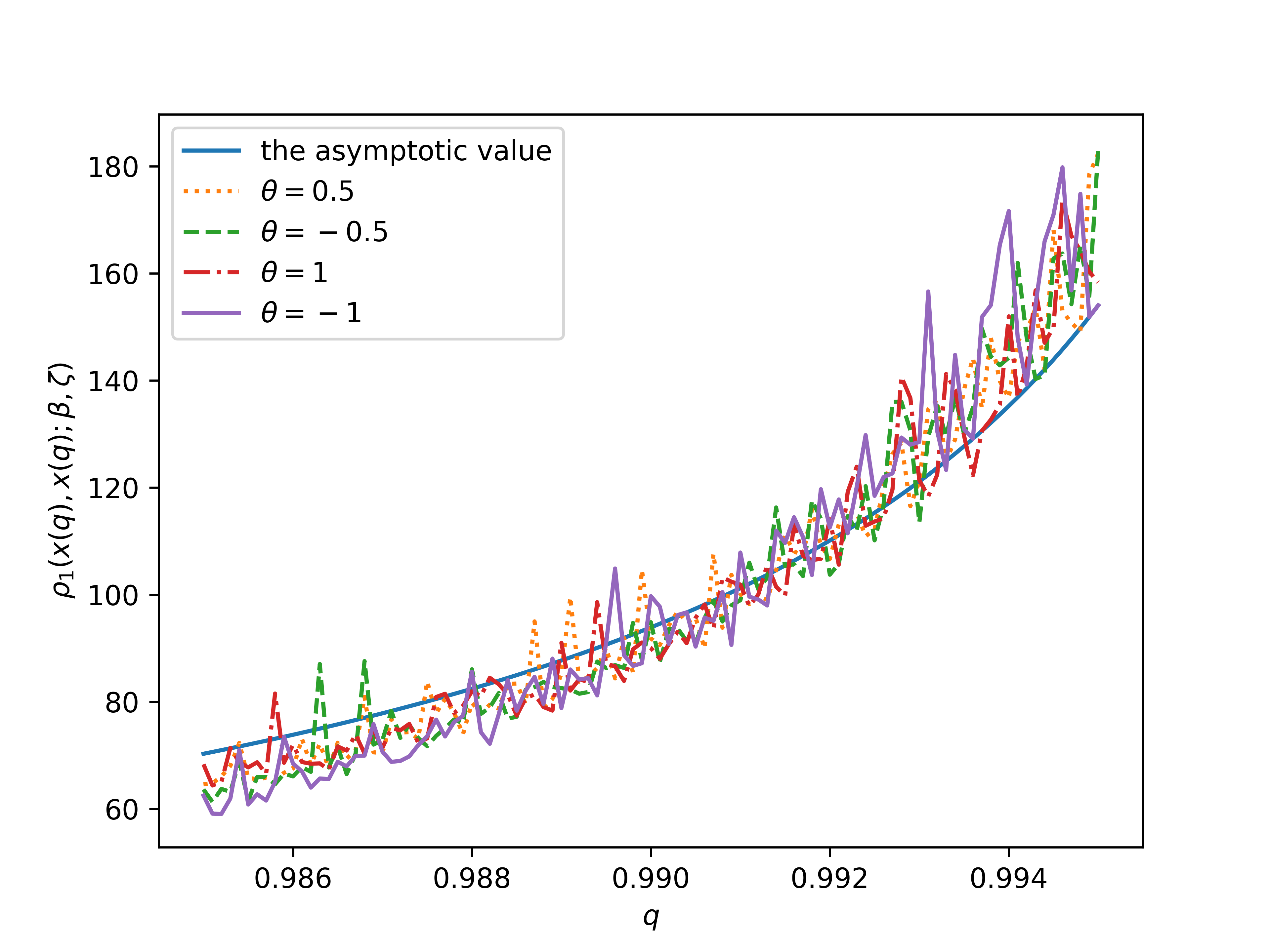} %插入图片，[]中设置图片大小，{}中是图片文件名
\captionsetup{font={footnotesize}}
\caption{Comparison between the asymptotic values for $\rho_1(x,x;\beta, \zeta)$  and its empirical values with $\beta=2$, $\zeta=\frac1{3}$, and $\alpha=2.8$.} % 最终文档中希望显示的图片标题
\label{Fig.main8} %用于文内引用的标签heavy_tailed_thtea
\end{figure}

\begin{figure}[htbp] %H为当前位置，!htb为忽略美学标准，htbp为浮动图形
\centering %图片居中
\includegraphics[width=0.8\textwidth]{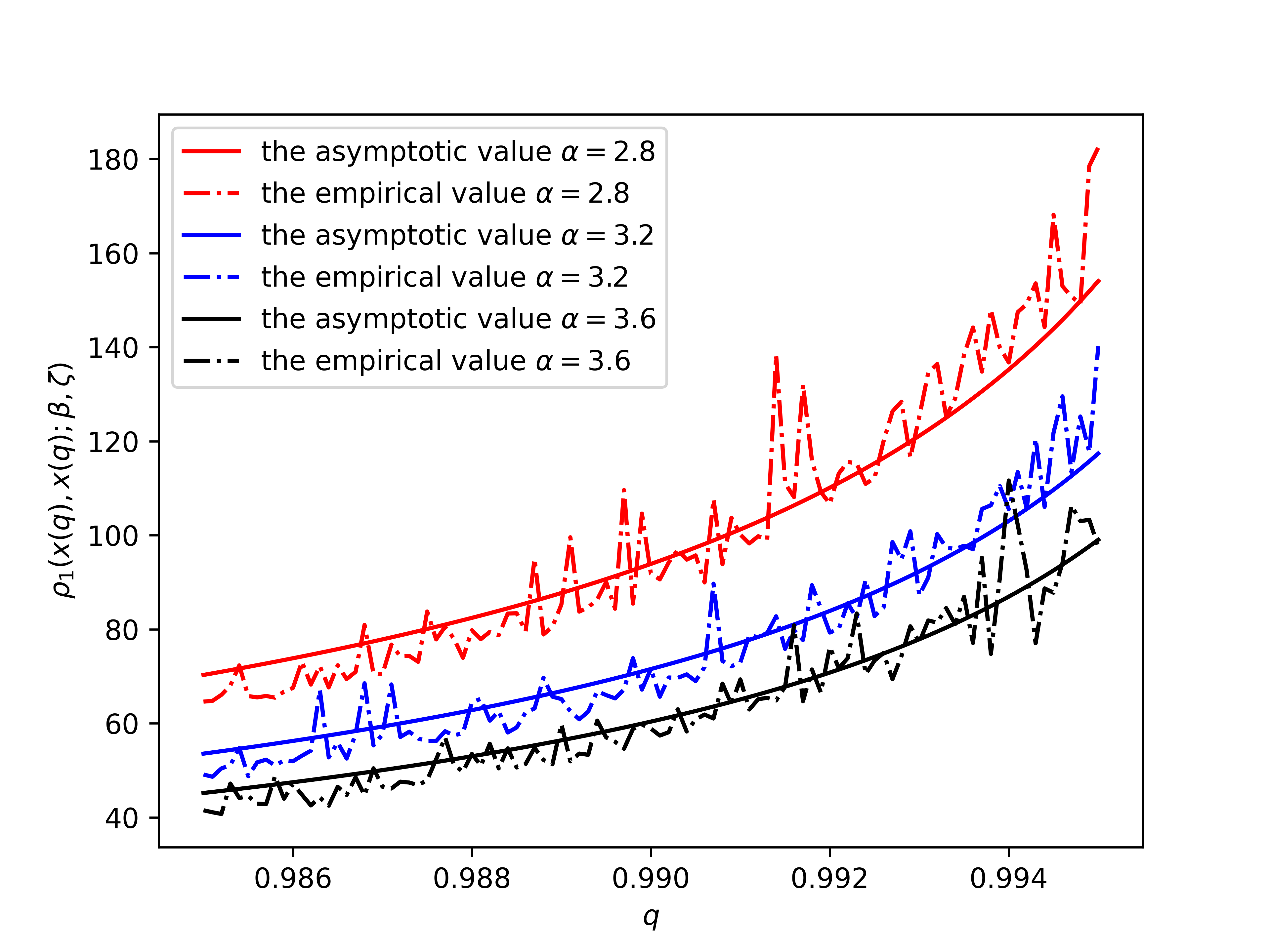} %插入图片，[]中设置图片大小，{}中是图片文件名
\captionsetup{font={footnotesize}}
\caption{Comparison between the asymptotic values for $\rho_1(x,x;\beta, \zeta)$  and its empirical values with $\beta=2$, $\zeta=\frac1{3}$, and $\theta=0.5$.} % 最终文档中希望显示的图片标题
\label{Fig.main9} %用于文内引用的标签heavy_tailed_thtea
\end{figure}

Consistent with Figures \ref{Fig.main4}-\ref{Fig.main6}, we assume that $\zeta=0$ in Tables \ref{tab1}-\ref{tab3}, and corresponding simulation results are summarized. In Figures \ref{Fig.main7}-\ref{Fig.main9}, we assume that $\zeta=\frac1{3}$. The results of the numerical simulation studies are presented in the form of line charts.

Some of the phenomena observed in Figures \ref{Fig.main4}-\ref{Fig.main6} are also evident in Figures \ref{Fig.main7}-\ref{Fig.main9}. Specifically, Figure  \ref{Fig.main7} highlights the amplification of errors caused by variations in $\beta$; Figure  \ref{Fig.main8} illustrates the influence of $\theta$ on estimation accuracy; and Figure  \ref{Fig.main9} demonstrates that an increase in $\alpha$ leads to larger errors, which aligns with the asymptotic formula (\ref{thm2res}).

%\subsection{Sensitivity analysis}

%\begin{table}[H]
%\centering
%\caption{Sensitivity analysis of asymptotic approximations.}
%\label{tabs}
%\resizebox{0.8\linewidth}{!}{
%\begin{tabular}{ccccc}
%\toprule[2pt]
%\multicolumn{2}{c}{\multirow{2}{*}{Model parameter $\zeta$ and its change}} & \multicolumn{3}{c}{asymptotic approximations} \\ \cline{3-5} 
%\multicolumn{2}{c}{}& $q=98.5\%$ & $q=99.0\%$ & $q=99.5\%$  \\ \hline
%\multirow{5}{*}{$\rho_1(x,y;\beta, \zeta)$, $\gamma=\frac{1}{2}$} 
%& $+2\%$   & 1   & 1     & 1     \\
%& $+1\%$   & 1   & 1     & 1     \\
%& ($\zeta=\frac1{3}$) & 1     & 1     & 1     \\
%& $-1\%$   & 1     & 1     & 1     \\
%& $-2\%$   & 1     & 1     & 1     \\ \hline
%\multirow{5}{*}{$\rho_1(x,y;\beta, \zeta)$, $\alpha=2.8$} 
%& $+2\%$   & 1     & 1     & 1     \\
%& $+1\%$   & 1     & 1     & 1     \\
%& ($\zeta=\frac1{3}$) & 1     & 1     & 1     \\
%& $-1\%$   & 1     & 1     & 1     \\
%& $-2\%$   & 1     & 1     & 1     \\ \bottomrule[2pt]
%\end{tabular}}
%\end{table}

In conclusion, Tables \ref{tab-1}-\ref{tab3} and Figures \ref{Fig.main4}-\ref{Fig.main9} demonstrate that as $q$ approaches 1 (or equivalently, as $x$ becomes relatively large), the asymptotic estimates of $\rho_1(x,y;\beta, \zeta)$ closely align with the computationally intensive simulated values. This indicates that our asymptotic results serve as accurate approximations and provide reliable insights.
%In conclusion, these numerical simulations demonstrate that our asymptotic formula performs well in estimating the real value of JTM. As one can see, when $x$ becomes relatively large, the more fluctuation the empirical value $\widehat{JTM}_1^{\beta}(x,x)$ exhibits. This may be due to the poor performance of the Monte Carlo method and we may repeat the simulation with a larger sample size $N$ to improve. We will not repeat it because it is beyond the scope of this study. Furthermore, the empirical value is overall larger than the asymptotic value in the case of heavy-tailedness when $x$ is relatively large and the case of light-tailedness is just the opposite. This reminds us that as the crisis approaches the obtained asymptotic formula may underestimate the risks when the risks are heavy-tailed and may overestimate them when they are light-tailed. However, our asymptotic result behaves excellently.  
\section{Proofs of the main results}\setcounter{thm}{0}\setcounter{equation}{0}\setcounter{lemma}{0}\setcounter{co}{0}\setcounter{re}{0}
Before we present lemmas and their proofs, we first introduce some notations that will facilitate our proofs. Let $\left\{X^{*}, X_i^{*},i\geq 1\right\}$ be a sequence of i.i.d. r.v.s with common distribution $F$,  $\left\{Y^{*}, Y_i^{*},i\geq 1\right\}$ be a sequence of i.i.d. r.v.s with common distribution $G$,  $\left\{X_{\wedge}^{*}, X_{\wedge,i}^{*},i\geq 1\right\}$ be a sequence of i.i.d. r.v.s with common distribution $F_{\wedge}$ satisfying $\overline{F_{\wedge}}(x)=\overline{F}^2(x)$ for all $x\in(-\infty,\infty)$,  and $\left\{Y_{\wedge}^{*}, Y_{\wedge,i}^{*},i\geq 1\right\}$ be a sequence of i.i.d. r.v.s with common distribution $G_{\wedge}$ satisfying $\overline{G_{\wedge}}(x)=\overline{G}^2(x)$ for all $x\in(-\infty,\infty)$. Furthermore, let the four sequences mentioned above be mutually independent and independent of any other randomness. With any fixed $n\geq 1$, we write $S_n^{*}=\sum_{i=1}^nX_i^{*}$, $T_n^{*}=\sum_{i=1}^nY_i^{*}$, $S_{\wedge,n}^{*}=\sum_{i=1}^nX_{\wedge, i}^{*}$, and $T_{\wedge, n}^{*}=\sum_{i=1}^nY_{\wedge, i}^{*}$. 
\subsection{Some lemmas}
The first lemma comes from Lemma 5.2 in {Pakes (2004)}, which reveals the asymptotic behavior of $n$-fold ($n\geq 2$) convolution for convolution-equivalent distributions.
\begin{lemma}\label{lemma-2}
	Let $F\in\mathscr{S}(\gamma)$, $\gamma\geq 0$, then for $n\geq 1$,
	\begin{align*}
		\lim\frac{\overline{F^{*n}}(x)}{\overline{F}(x)}= n\left(\widehat{F}(\gamma)\right)^{n-1}.
	\end{align*}
\end{lemma}
The establishment of the following lemma is initially motivated by Lemma 3.2 in {Tang and Tsitsiashvili (2003)}. It shows the closure property under tail equivalence for the convolution-equivalent distribution class.
\begin{lemma}\label{lemma-1}
	Let $F$, $F_1$, and $F_2$ be three distributions on $(-\infty,\infty)$ such that $F\in\mathscr{S}(\gamma)$, $\gamma\geq 0$, and that the limit $l_i=\lim\frac{\overline{F_i}(x)}{\overline{F}(x)}$ exists for $i=1,2$. Then
	\begin{align}\label{lemma-1con}
		\lim\frac{\overline{F_1*F_2}(x)}{\overline{F}(x)}=l_1\widehat{F_2}(\gamma)+l_2\widehat{F_1}(\gamma).
	\end{align}
\end{lemma}
{\bf Proof}\quad The proof primarily involves dividing the argument into three distinct cases. {\bf Case 1: $l_1\neq 0$ and $l_2\neq 0$}. In this scenario, $F_1$, $F_2\in\mathscr{S}(\gamma)$ due to Theorem 2.7 in Embrechts and Goldie (1982) and Corollary 2.1 (i) in Pakes (2004), then the desired conclusion is an immediate consequence of Lemma 3.2 in Tang and Tsitsiashvili (2003). 

\noindent {\bf Case 2: $l_1=l_2=0$}. For any fixed $0<\varepsilon<1$, there exists some $x_0>0$ such that for any $x>x_0$, $\overline{F_i}(x)\leq \varepsilon\overline{F}(x)$ holds for $i=1,2$. We introduce a survival distribution defined as 
	\begin{align}\label{lemma0prove1}
		\overline{F_{\varepsilon}}(x)=\max\left\{\overline{F_1}(x),\overline{F_2}(x),\varepsilon\overline{F}(x)\right\}
	\end{align}
 and it is obvious that $F_{\varepsilon}\in\mathscr{S}(\gamma)$ due to $\overline{F_{\varepsilon}}(x)\sim \varepsilon\overline{F}(x)$.
% and there exists some positive constant $C$, which is independent of $\varepsilon$ and $x$, such that $\overline{F_{\varepsilon}}(x)\leq C\overline{F}(x)$ for all $x\geq 0$. T
Then according to (\ref{lemma0prove1}) and the result derived in {\bf Case 1}, we derive that
	\begin{align*}
		0\leq \lim\frac{\overline{F_1*F_2}(x)}{\overline{F}(x)}&\leq \lim\frac{\overline{F_{\varepsilon}*F_{\varepsilon}}(x)}{\overline{F}(x)}=2\varepsilon\widehat{F_{\varepsilon}}(\gamma)\\
		&\leq 2\varepsilon\left(\widehat{F_1}(\gamma)+\widehat{F_2}(\gamma)+\varepsilon\widehat{F}(\gamma)\right)\\
		&\to 0, \quad \rm{as} \quad \varepsilon\downarrow 0,
	\end{align*}
	which meets (\ref{lemma-1con}). 
	
\noindent	{\bf Case 3: $l_1=0$, $l_2\neq 0$ or $l_1\neq 0$, $l_2=0$}. Without loss of generality, let $l_1\neq 0$, $l_2=0$. Following the methodology utilized in {\bf Case 2}, where we introduced a survival distribution, the asymptotic upper bound is readily established. Regarding the asymptotic lower bound, we have for some $K>0$ that
	\begin{align*}
		\overline{F_1*F_2}(x)&\geq \int_{-K}^{K}\overline{F_1}(x-y)F_2(\mathrm{d}y)\\
		&\sim l_1\int_{-K}^{K}e^{\gamma y}F_2(\mathrm{d}y)\overline{F}(x)\\
		&\sim l_1\widehat{F_2}(\gamma)\overline{F}(x), \quad {\rm as} \quad K\to\infty.
	\end{align*}
	 This ends the proof of Lemma \ref{lemma-1}.$\hfill\blacksquare$ 

The following lemma can be viewed as a corollary of Lemma \ref{lemma-1}. Nevertheless, we present it as an independent lemma due to its frequent application in the following proofs. The lemma roughly conveys that if a convolution-equivalent distribution dominates a distribution, then the $n$-fold convolution of this distribution is still dominated by the same convolution-equivalent distribution.
\begin{lemma}\label{lemma0}
	Let $\left\{X_i,i\geq 1\right\}$ be a sequence of i.i.d. r.v.s with common distribution $G$. If there exists some distribution $F\in\mathscr{S}(\gamma)$, $\gamma\geq 0$, satisfying $\overline{G}(x)=o\left(\overline{F}(x)\right)$, then for any fixed $n\geq 1$, it holds that $P\left(\sum_{i=1}^n X_i>x\right)=o\left(\overline{F}(x)\right)$.
\end{lemma}

The next lemma provides an asymptotic estimate of the joint tail probability of $S_n$ and $T_n$ with $n$ being a positive integer, which also implies that $S_n$ and $T_n$ are strongly asymptotically independent. Before presenting the lemma, let us introduce some notations. For $\gamma\geq 0$, $\theta\in[-1,1]$, $n_1,n_2,n_3,n_4,n\in\mathbb{N}$, and $n_1+n_2+n_3+n_4=n$, we write
\begin{align*}
		\left\{\begin{aligned}
	%&\Omega(\gamma,\theta)=(1+\theta)\widehat{F}(\gamma)\widehat{G}(\gamma)-\theta \widehat{F_{\wedge}}(\gamma)\widehat{G}(\gamma)-\theta\widehat{F}(\gamma)\widehat{G_{\wedge}}(\gamma)+\theta \widehat{F_{\wedge}}(\gamma)\widehat{G_{\wedge}}(\gamma)\\
	&A^{\theta}_{n_1,n_2,n_3,n_4}\left(n\right)=\frac{n!}{n_1!n_2!n_3!n_4!}(-1)^{n_2+n_3}(1+\theta)^{n_1}\theta^{n-n_1}\\
	&B^{\gamma,F}_{n_1,n_2,n_3,n_4}\left(n\right)=(n_1+n_3)\left(\widehat{F}(\gamma)\right)^{n_1+n_3-1}\left(\widehat{F_{\wedge}}(\gamma)\right)^{n_2+n_4}\\
	&\Gamma^{\gamma,G}_{n_1,n_2,n_3,n_4}\left(n\right)=(n_1+n_2)\left(\widehat{G}(\gamma)\right)^{n_1+n_2-1}\left(\widehat{G_{\wedge}}(\gamma)\right)^{n_3+n_4}\\
	&K(n,\theta,\gamma)= \frac{\partial^2f^n(u,v,s,t)}{\partial u\partial v}\Bigg\vert_{u=\widehat{F}(\gamma),v=\widehat{G}(\gamma),s=\widehat{F_{\wedge}}(\gamma),t=\widehat{G_{\wedge}}(\gamma)},
\end{aligned}\right.
\end{align*}
where $f(u,v,s,t)=(1+\theta)uv-\theta sv-\theta ut+\theta st$.
\begin{lemma}\label{lemma1}
	Let $\left\{(X,Y), (X_i, Y_i),i\geq 1\right\}$ be a sequence of i.i.d. random pairs with marginal distributions $F, G\in\mathscr{S}(\gamma)$, $\gamma\geq 0$, respectively. Let the joint distribution of $(X, Y)$ satisfy (\ref{joint}) with parameter $\theta\in [-1,1]$, then we have 
	\begin{align}\label{lemma1res}
		P\left(S_n>x, T_n>y\right)\sim K(n,\theta,\gamma)\overline{F}(x)\overline{G}(y)
	\end{align}
	holds for any fixed $n\geq 1$, where $f\sim 0\cdot g$ is understood as $f=o(g)$ for two positive functions $f$, $g$.
\end{lemma}
\begin{re}\label{re4.1}
With some obvious calculations, we get $K(1,-1,0)=0$, then $$P\left(S_n>x, T_n>y\right)=o\left(\overline{F}(x)\overline{G}(y)\right)$$ when $n=1$ and $\theta=-1$, then Lemma \ref{lemma1} fails to give a precise asymptotic estimate of $P\left(S_n>x, T_n>y\right)$. While we claim in this remark that the coefficient $K(n,\theta,\gamma)>0$ for all $n\geq 2$ and $\theta\in[-1,1]$. It is obvious that with $n\geq 2$,  $K(n,\theta,\gamma)>0$ holds for $\gamma=0$, and in the following we consider the case of $\gamma>0$. By direct computation, for $\gamma>0$ and non-degenerate $F$ and $G$, we have the following iterative formula.
	 \begin{align*}
	 	K(n,\theta,\gamma)&=(1+\theta)nf^{n-1}\left(\widehat{F}(\gamma),\widehat{G}(\gamma),\widehat{F_{\wedge}}(\gamma),\widehat{G_{\wedge}}(\gamma)\right)\\
	 	&\quad +(n-1)n\left(\widehat{F}(\gamma)+\theta\left(\widehat{F}(\gamma)-\widehat{F_{\wedge}}(\gamma)\right)\right)\\
	 	&\quad\times \left(\widehat{G}(\gamma)+\theta\left(\widehat{G}(\gamma)-\widehat{G_{\wedge}}(\gamma)\right)\right)f^{n-2}\left(\widehat{F}(\gamma),\widehat{G}(\gamma),\widehat{F_{\wedge}}(\gamma),\widehat{G_{\wedge}}(\gamma)\right).
	 \end{align*}
	By noticing that $\widehat{F}(\gamma)>\widehat{F_{\wedge}}(\gamma)>0$, $\widehat{G}(\gamma)> \widehat{G_{\wedge}}(\gamma)>0$, and 
	 \begin{align*}
	 	f\left(\widehat{F}(\gamma),\widehat{G}(\gamma),\widehat{F_{\wedge}}(\gamma),\widehat{G_{\wedge}}(\gamma)\right)> \left(\widehat{F}(\gamma)+\theta\left(\widehat{F}(\gamma)-\widehat{F_{\wedge}}(\gamma)\right)\right)\left(\widehat{G}(\gamma)-\widehat{G_{\wedge}}(\gamma)\right).
	 \end{align*} 
	 To prove $K(n,\theta,\gamma)>0$, we only need to prove 
	 \begin{align}\label{rmpr}
	 	\widehat{H}(\gamma)+\theta\left(\widehat{H}(\gamma)-\widehat{H_{\wedge}}(\gamma)\right)>0
	 \end{align}
	 with $H$ to be $F$ or $G$. It is clear that (\ref{rmpr}) holds for $\theta\in[0,1]$. While in the case of $\theta\in[-1,0)$, it is also evident by noticing that $\widehat{H}(\gamma)+\theta\left(\widehat{H}(\gamma)-\widehat{H_{\wedge}}(\gamma)\right)=(1+\theta)\widehat{H}(\gamma)-\theta\widehat{H_{\wedge}}(\gamma)$. This completes the remark.
\end{re}
{\bf Proof}\quad 
We use mathematical induction to prove the following statement. For any $n\in \mathbb{N}$, it holds for all $x,y\in(-\infty,\infty)$ that 
\begin{align}\label{lemma1re2}
	P(S_n>x, T_n>y)&=\sum_{\substack{n_1,n_2,n_3,n_4\in\mathbb{N},\\n_1+n_2+n_3+n_4=n}}A^{\theta}_{n_1,n_2,n_3,n_4}\left(n\right)\nonumber\\
	&\quad\times P\left(S_{n_1+n_3}^{*}+S_{\wedge,n_2+n_4}^{*}>x\right)P\left(T_{n_1+n_2}^{*}+T_{\wedge,n_3+n_4}^{*}>y\right),
\end{align}
which is an extension of relation (\ref{decomp}). It is not hard to verify by (\ref{decomp}) that the statement holds for $n=1$. Assuming that the statement is true for some $n$, we need to prove that the statement holds for $n+1$. By conditioning on $(S_n, T_n)$ and using (\ref{decomp}), we derive
%\allowdisplaybreaks
\begin{align}\label{delta0}
	&\quad P(S_{n+1}>x, T_{n+1}>y)\nonumber\\
	&=(1+\theta)P\left(S_n+X_{n+1}^{*}>x,T_n+Y_{n+1}^{*}>y\right)-\theta P\left(S_n+X_{\wedge,n+1}^{*}>x,T_n+Y_{n+1}^{*}>y\right)\nonumber\\
	&\quad -\theta P\left(S_n+X_{n+1}^{*}>x,T_n+Y_{\wedge,n+1}^{*}>y\right)+\theta P\left(S_n+X_{\wedge,n+1}^{*}>x,T_n+Y_{\wedge,n+1}^{*}>y\right)\nonumber\\
	&\triangleq \Delta_1(x,y;\theta)+\Delta_2(x,y;\theta)+\Delta_3(x,y;\theta)+\Delta_4(x,y;\theta).
\end{align}
Utilizing the inductive hypothesis, we have
\begin{align}\label{delta1}
	\Delta_1(x,y;\theta)&=(1+\theta)\sum_{\substack{n_1,n_2,n_3,n_4\in\mathbb{N},\\n_1+n_2+n_3+n_4=n}}A^{\theta}_{n_1,n_2,n_3,n_4}\left(n\right)\nonumber\\
	&\quad \times P\left(S_{n_1+1+n_3}^{*}+S_{\wedge,n_2+n_4}^{*}>x\right)P\left(T_{n_1+1+n_2}^{*}+T_{\wedge,n_3+n_4}^{*}>y\right)\nonumber\\
	&=\sum_{\substack{n_1,n_2,n_3,n_4\in\mathbb{N},\\n_1+n_2+n_3+n_4=n+1}}\frac{n!}{(n_1-1)!n_2!n_3!n_4!}(-1)^{n_2+n_3}(1+\theta)^{n_1}\theta^{(n+1)-n_1}\nonumber\\
	&\quad \times P\left(S_{n_1+n_3}^{*}+S_{\wedge,n_2+n_4}^{*}>x\right)P\left(T_{n_1+n_2}^{*}+T_{\wedge,n_3+n_4}^{*}>y\right).
\end{align}
By analogous methods, we can also derive
\begin{align}\label{delta2}
	\Delta_2(x,y;\theta)&=\sum_{\substack{n_1,n_2,n_3,n_4\in\mathbb{N},\\n_1+n_2+n_3+n_4=n+1}}\frac{n!}{n_1!(n_2-1)!n_3!n_4!}(-1)^{n_2+n_3}(1+\theta)^{n_1}\theta^{(n+1)-n_1}\nonumber\\
	&\quad \times P\left(S_{n_1+n_3}^{*}+S_{\wedge,n_2+n_4}^{*}>x\right)P\left(T_{n_1+n_2}^{*}+T_{\wedge,n_3+n_4}^{*}>y\right),
\end{align}
\begin{align}\label{delta3}
	\Delta_3(x,y;\theta)&=\sum_{\substack{n_1,n_2,n_3,n_4\in\mathbb{N},\\n_1+n_2+n_3+n_4=n+1}}\frac{n!}{n_1!n_2!(n_3-1)!n_4!}(-1)^{n_2+n_3}(1+\theta)^{n_1}\theta^{(n+1)-n_1}\nonumber\\
	&\quad \times P\left(S_{n_1+n_3}^{*}+S_{\wedge,n_2+n_4}^{*}>x\right)P\left(T_{n_1+n_2}^{*}+T_{\wedge,n_3+n_4}^{*}>y\right),
\end{align}
and 
\begin{align}\label{delta4}
	\Delta_4(x,y;\theta)&=\sum_{\substack{n_1,n_2,n_3,n_4\in\mathbb{N},\\n_1+n_2+n_3+n_4=n+1}}\frac{n!}{n_1!n_2!n_3!(n_4-1)!}(-1)^{n_2+n_3}(1+\theta)^{n_1}\theta^{(n+1)-n_1}\nonumber\\
	&\quad \times P\left(S_{n_1+n_3}^{*}+S_{\wedge,n_2+n_4}^{*}>x\right)P\left(T_{n_1+n_2}^{*}+T_{\wedge,n_3+n_4}^{*}>y\right).
\end{align}
By substituting (\ref{delta1})-(\ref{delta4}) into (\ref{delta0}) and utilizing  the elementary equation $$\frac{(n+1)!}{n_1!n_2!n_3!n_4!}=\frac{n!}{(n_1-1)!n_2!n_3!n_4!}+\frac{n!}{n_1!(n_2-1)!n_3!n_4!}+\frac{n!}{n_1!n_2!(n_3-1)!n_4!}+\frac{n!}{n_1!n_2!n_3!(n_4-1)!},$$ we establish 
\begin{align*}
	&\quad P(S_{n+1}>x, T_{n+1}>y)\nonumber\\
	&=\sum_{\substack{n_1,n_2,n_3,n_4\in\mathbb{N},\\n_1+n_2+n_3+n_4=n+1}}A^{\theta}_{n_1,n_2,n_3,n_4}\left(n+1\right)P\left(S_{n_1+n_3}^{*}+S_{\wedge,n_2+n_4}^{*}>x\right)P\left(T_{n_1+n_2}^{*}+T_{\wedge,n_3+n_4}^{*}>y\right),
\end{align*} 
which allows us to conclude that (\ref{lemma1re2}) holds for any $n\in \mathbb{N}$. After establishing the previous statement, for $n_1$, $n_2$, $n_3$, and $n_4$ indicated in (\ref{lemma1re2}), by applying Lemmas  \ref{lemma-2}-\ref{lemma0}, we deduce that 
\begin{align}\label{delta5}
	P\left(S_{n_1+n_3}^{*}+S_{\wedge,n_2+n_4}^{*}>x\right)\sim (n_1+n_3)\left(\widehat{F}(\gamma)\right)^{n_1+n_3-1}\left(\widehat{F_{\wedge}}(\gamma)\right)^{n_2+n_4}\overline{F}(x)
\end{align}
and
\begin{align}\label{delta6}
	P\left(T_{n_1+n_2}^{*}+T_{\wedge,n_3+n_4}^{*}>y\right)\sim (n_1+n_2)\left(\widehat{G}(\gamma)\right)^{n_1+n_2-1}\left(\widehat{G_{\wedge}}(\gamma)\right)^{n_3+n_4}\overline{G}(y).
\end{align}
Substituting (\ref{delta5}) and (\ref{delta6}) into (\ref{lemma1re2}) gives 
\begin{align*}
	P(S_{n}>x, T_{n}>y)&\sim \sum_{\substack{n_1,n_2,n_3,n_4\in\mathbb{N},\\n_1+n_2+n_3+n_4=n}}A^{\theta}_{n_1,n_2,n_3,n_4}\left(n\right)B^{\gamma,F}_{n_1,n_2,n_3,n_4}\left(n\right)\Gamma^{\gamma,G}_{n_1,n_2,n_3,n_4}\left(n\right) \overline{F}(x)\overline{G}(y)\\
	&=K(n,\theta,\gamma)\overline{F}(x)\overline{G}(y),
\end{align*}
where the last step is due to the following elementary equation
\begin{align*}
	f^n(u,v,s,t)=\sum_{\substack{n_1,n_2,n_3,n_4\in\mathbb{N},\\n_1+n_2+n_3+n_4=n}}A^{\theta}_{n_1,n_2,n_3,n_4}\left(n\right)u^{n_1+n_3}v^{n_1+n_2}s^{n_2+n_4}t^{n_3+n_4}.
\end{align*}
This ends the proof of Lemma \ref{lemma1}.
$\hfill\blacksquare$ 

\begin{re}\label{re4.2}
Under the conditions of Lemma \ref{lemma1}, this remark will show 
\begin{align}\label{re4.22}
	P\left(S_{n}>x, T_{m}>y\right)\asymp\overline{F}(x)\overline{G}(y)
\end{align}
for $n\geq 1$, $m\geq 2$ or $n\geq 2$, $m\geq 1$. Without loss of generality, we assume $n\geq 1$, $m\geq 2$, and $m\geq n$. Then, for $n,m\geq 2$, by Lemma \ref{lemma1} and Remark \ref{re4.1}, there exists a positive constant $C>1$ such that 
\begin{align}\label{re4.21}
	\frac{1}{C}\overline{F}(x)\overline{G}(y)\leq P\left(S_{n}>x, T_{n}>y\right)\leq C\overline{F}(x)\overline{G}(y)
\end{align}
holds for all $x,y\in(-\infty, \infty)$. It can be deduced by (\ref{re4.21}) and $G\in\mathscr{S}(\gamma)$ that
\begin{align*}
	P\left(S_{n}>x, T_{m}>y\right)\leq C\overline{F}(x)P(T_{m-n+1}>y)=O\left(\overline{F}(x)\overline{G}(y)\right).
\end{align*}
Following the analogous argument, one can derive that (\ref{re4.22}) holds for $n,m\geq 2$. As for $n=1$ and $m\geq 2$, according to (\ref{decomp}), we have
\begin{align*}
	P(X_1>x, T_n>y)&=(1+\theta)P(X^{*}>x, T_{n-1}+Y^{*}>y)-\theta P(X^{*}>x, T_{n-1}+Y^{*}_{\wedge}>y)\\
	&\quad -\theta P(X^{*}_{\wedge}>x, T_{n-1}+Y^{*}>y)+\theta P(X^{*}_{\wedge}>x, T_{n-1}+Y^{*}_{\wedge}>y)\\
	&\sim \left((1+\theta)n\left(\widehat{G}(\gamma)\right)^{n-1}-\theta(n-1)\left(\widehat{G}(\gamma)\right)^{n-2}\widehat{G_{\wedge}}(\gamma)\right)\overline{F}(x)\overline{G}(y),
\end{align*}
which indicates (\ref{re4.22}) by noticing $(1+\theta)n\left(\widehat{G}(\gamma)\right)^{n-1}-\theta(n-1)\left(\widehat{G}(\gamma)\right)^{n-2}\widehat{G_{\wedge}}(\gamma)>0$. While an exact asymptotic equivalence for $P\left(S_{n}>x, T_{m}>y\right)$ can be established in a manner analogous to Lemma \ref{lemma1}, we omit its explicit form due to the cumbersome notation involved. Nevertheless, (\ref{re4.22}) established in this remark suffices to prove the main results.
\end{re}
The following lemma further explores the joint tail probability of $S_n$ and $T_n$ when one of $X_i$ is large, thereby laying the foundation for establishing the main results.
\begin{lemma}\label{lemma2}
Let the conditions of Lemma \ref{lemma1} be valid. Then for any $0<u<1$ and $n\geq 2$, it holds for $1\leq i \leq n$ that
\begin{align}\label{lemma2res}
P\left(X_i>u x,S_n>x, T_n>y\right)\sim \frac{1}{n}P\left(S_n>x, T_n>y\right).
\end{align}
Particularly, if $F$, $G\in\mathscr{R}_{-\alpha}$ for some $\alpha\geq 0$, then it holds uniformly for $u\geq\varepsilon$ with some $\varepsilon\in(0,1)$ that
 \begin{align}\label{lemma2res2}
P\left(X_i>ux,S_n>x, T_n>y\right)\sim (n+\theta)\overline{F}((u\vee 1)x).\overline{G}(y)
\end{align}
If further $F(-x)=O\left(\overline{F}(x)\right)$, then we have
 \begin{align}\label{lemma2res3}
 	P\left(X_i^{-}>ux,S_n>x, T_n>y\right)=o(1)\overline{F}((u\vee 1)x)\overline{G}(y)
 \end{align}
 	uniformly for $u\geq\varepsilon$ with some $\varepsilon\in(0,1)$.
\end{lemma}
{\bf Proof}\quad
%It is sufficient to prove (\ref{lemma2res}) because (\ref{lemma2res1}) is obtained by adopting Lemma \ref{lemma1}. 
Without loss of generality, let $i=1$. Given (\ref{decomp}), for fixed $n\geq 2$ and $0<u<1$, we do the following decomposition.
\begin{align}\label{lemma2re1}
	P\left(X_1>ux,S_n>x, T_n>y\right)&=(1+\theta)P\left(X_1^{*}>ux,S_{n-1}+X_1^{*}>x, T_{n-1}+Y_1^{*}>y\right)\nonumber\\
	&\quad -\theta P\left(X_1^{*}>ux,S_{n-1}+X_1^{*}>x, T_{n-1}+Y_{\wedge,1}^{*}>y\right)\nonumber\\
	&\quad -\theta P\left(X_{\wedge,1}^{*}>ux,S_{n-1}+X_{\wedge,1}^{*}>x, T_{n-1}+Y_1^{*}>y\right)\nonumber\\
	&\quad +\theta P\left(X_{\wedge,1}^{*}>ux,S_{n-1}+X_{\wedge,1}^{*}>x, T_{n-1}+Y_{\wedge,1}^{*}>y\right)\nonumber\\
	&\triangleq I_1(x,y)+ I_2(x,y)+I_3(x,y)+I_4(x,y).
\end{align}
By (\ref{lemma1re2}), we obtain 
\begin{align*}
	I_1(x,y)&=\sum_{\substack{n_1,n_2,n_3,n_4\in\mathbb{N},\\n_1+n_2+n_3+n_4=n-1}}\frac{(n-1)!}{n_1!n_2!n_3!n_4!}(1+\theta)^{n_1+1}(-\theta)^{n_2+n_3}\theta^{n_4}\nonumber\\
		&\quad\times P\left(X_1^{*}>ux, S_{n_1+1+n_3}^{*}+S_{\wedge,n_2+n_4}^{*}>x\right)P\left(T_{n_1+1+n_2}^{*}+T_{\wedge,n_3+n_4}^{*}>y\right).
\end{align*}
It follows from (\ref{delta6}) that
\begin{align}\label{lemma2re3}
	P\left(T_{n_1+1+n_2}^{*}+T_{\wedge,n_3+n_4}^{*}>y\right)\sim (n_1+n_2+1)\left(\widehat{G}(\gamma)\right)^{n_1+n_2}\left(\widehat{G_{\wedge}}(\gamma)\right)^{n_3+n_4}\overline{G}(y).
\end{align}
%where $Z^{*}$ can be $X_{1}^{*}$ or $X_{\wedge,1}^{*}$. If $Z^{*}=X_{1}^{*}$, then for every fixed $c>0$, it holds for large $x$ that
For every fixed $c>0$, it holds for large $x$ that
\begin{align}\label{lemma2re4}
	&\quad P\left(X_1^{*}>ux, S_{n_1+1+n_3}^{*}+S_{\wedge,n_2+n_4}^{*}>x\right)\nonumber\\
	&=P\left(S_{n_1+1+n_3}^{*}+S_{\wedge,n_2+n_4}^{*}>x\right)-P\left(X_1^{*}\leq c, S_{n_1+1+n_3}^{*}+S_{\wedge,n_2+n_4}^{*}>x\right)\nonumber\\
	&\quad -P\left(c<X_1^{*}\leq ux, S_{n_1+1+n_3}^{*}+S_{\wedge,n_2+n_4}^{*}>x\right)\nonumber\\
	&\triangleq P_1(x)-P_2(x)-P_3(x).
\end{align}
Utilizing (\ref{delta5}), it holds that
\begin{align}\label{lemma2p1}
	P_1(x)\sim (n_1+1+n_3)\left(\widehat{F}(\gamma)\right)^{n_1+n_3}\left(\widehat{F_{\wedge}}(\gamma)\right)^{n_2+n_4}\overline{F}(x).
\end{align}
By employing (\ref{delta5}), we obtain
\begin{align}\label{lemma2p2}
	P_2(x)&=\int_{-\infty}^cP\left(S_{n_1+n_3}^{*}+S_{\wedge,n_2+n_4}^{*}>x-u\right)F(\mathrm{d}u)\nonumber\\
	&\sim (n_1+n_3)\left(\widehat{F}(\gamma)\right)^{n_1+n_3-1}\left(\widehat{F_{\wedge}}(\gamma)\right)^{n_2+n_4}\int_{-\infty}^c\overline{F}(x-u)F(\mathrm{d}u)\nonumber\\
	&\sim (n_1+n_3)\left(\widehat{F}(\gamma)\right)^{n_1+n_3-1}\left(\widehat{F_{\wedge}}(\gamma)\right)^{n_2+n_4}\overline{F}(x)\int_{-\infty}^ce^{\gamma u}F(\mathrm{d}u)\nonumber\\
	&\sim (n_1+n_3)\left(\widehat{F}(\gamma)\right)^{n_1+n_3}\left(\widehat{F_{\wedge}}(\gamma)\right)^{n_2+n_4}\overline{F}(x), \quad {\rm as}\quad c\to\infty,
\end{align}
where, in the third step, we applied that $\lim\frac{\overline{F}(x-u)}{\overline{F}(x)}=e^{\gamma u}$ holds uniformly for $u\leq c$ if $\gamma>0$  and the dominated convergence theorem if $\gamma=0$. As for $P_3(x)$, by noticing that $P(X^{*}_{\wedge}>x)\leq P(X^{*}>x)$ holds for all $x\in (-\infty, \infty)$, we derive that
\begin{align}\label{lemma2p3}
	 \lim_{c\to\infty}\lim \frac{P_3(x)}{\overline{F}(x)}
	&\leq \lim_{c\to\infty}\lim\frac{P\left(c<X_1^{*}\leq ux, S_{n}^{*}>x\right)}{\overline{F}(x)}\nonumber\\
	&\leq\lim_{c\to\infty}\lim \frac{P\left(S_{n}^{*}>x, \sum_{k=2}^nX_{k}^{*}>(1-u)x, X_1^{*}>c\right)}{\overline{F}(x)}\nonumber\\
	&\leq \sum_{k=2}^{n}\lim_{c\to\infty}\lim\frac{P\left(S_{n}^{*}>x, X_k^{*}>c, X_1^{*}>c\right)}{\overline{F}(x)}\nonumber\\
	&=0,
\end{align}
where we used Lemma 3.1 in {Li and Tang (2015)} in the last step. Plugging (\ref{lemma2p1})-(\ref{lemma2p3}) into (\ref{lemma2re4}) yields that 
\begin{align}\label{lemma2re5}
	P\left(X_1^{*}>ux, S_{n_1+1+n_3}^{*}+S_{\wedge,n_2+n_4}^{*}>x\right)&\sim \left(\widehat{F}(\gamma)\right)^{n_1+n_3}\left(\widehat{F_{\wedge}}(\gamma)\right)^{n_2+n_4}\overline{F}(x).
\end{align}
A combination of (\ref{lemma2re3}) and (\ref{lemma2re5}) gives that
\begin{align}\label{lemma2re6}
	&\quad I_1(x,y)\nonumber\\
	 &\sim \sum_{\substack{n_1,n_2,n_3,n_4\in\mathbb{N},\\n_1+n_2+n_3+n_4=n-1}}\frac{(n-1)!}{n_1!n_2!n_3!n_4!}(1+\theta)^{n_1+1}(-\theta)^{n_2+n_3}\theta^{n_4}(n_1+n_2+1)\nonumber\\
	&\quad\times \left(\widehat{F}(\gamma)\right)^{n_1+n_3}\left(\widehat{F_{\wedge}}(\gamma)\right)^{n_2+n_4}\left(\widehat{G}(\gamma)\right)^{n_1+n_2}\left(\widehat{G_{\wedge}}(\gamma)\right)^{n_3+n_4}\overline{F}(x)\overline{G}(y)\nonumber\\
	&=\frac{1}{n}\sum_{\substack{n_1,n_2,n_3,n_4\in\mathbb{N},\\n_1+n_2+n_3+n_4=n}}\frac{n!}{n_1!n_2!n_3!n_4!}(1+\theta)^{n_1}(-\theta)^{n_2+n_3}\theta^{n_4}n_1(n_1+n_2)\nonumber\\
	&\quad\times \left(\widehat{F}(\gamma)\right)^{n_1+n_3-1}\left(\widehat{F_{\wedge}}(\gamma)\right)^{n_2+n_4}\left(\widehat{G}(\gamma)\right)^{n_1+n_2-1}\left(\widehat{G_{\wedge}}(\gamma)\right)^{n_3+n_4}\overline{F}(x)\overline{G}(y).
\end{align}
The same steps can be used as in the verification of (\ref{lemma2re6}) to obtain
\begin{align}\label{lemma2re7}
	&\quad I_2(x,y)\nonumber\\
	&\sim\frac1{n}\sum_{\substack{n_1,n_2,n_3,n_4\in\mathbb{N},\\n_1+n_2+n_3+n_4=n}}\frac{n!}{n_1!n_2!n_3!n_4!}(1+\theta)^{n_1}(-\theta)^{n_2+n_3}\theta^{n_4}n_3(n_1+n_2)\nonumber\\
	&\quad\times \left(\widehat{F}(\gamma)\right)^{n_1+n_3-1}\left(\widehat{F_{\wedge}}(\gamma)\right)^{n_2+n_4}\left(\widehat{G}(\gamma)\right)^{n_1+n_2-1}\left(\widehat{G_{\wedge}}(\gamma)\right)^{n_3+n_4}\overline{F}(x)\overline{G}(y),
\end{align}
\begin{align}\label{lemma2re8}
	I_3(x,y)=o\left(\overline{F}(x)\overline{G}(y)\right),
\end{align}
and 
\begin{align}\label{lemma2re9}
	I_4(x,y)=o\left(\overline{F}(x)\overline{G}(y)\right).
\end{align}
Plugging (\ref{lemma2re6})-(\ref{lemma2re9}) into (\ref{lemma2re1}) and employing Lemma \ref{lemma1} yield that
\begin{align*}
P\left(X_1>ux,S_n>x, T_n>y\right)\sim \frac1{n}K(n,\theta,\gamma)\overline{F}(x)\overline{G}(y)\sim \frac1{n}P\left(S_n>x, T_n>y\right).
\end{align*}
This ends the proof of (\ref{lemma2res}).

Next, we prove the relation (\ref{lemma2res2}). On the one hand, we have uniformly for $\varepsilon\leq u\leq 1$ that
\begin{align}\label{Szudd}
	P\left(X_1>ux,S_n>x, T_n>y\right)&\leq P\left(X_1>\varepsilon x,S_n>x, T_n>y\right)\nonumber\\
	&\sim \frac{1}{n}P\left(S_n>x, T_n>y\right)\nonumber\\
	&\sim (n+\theta)\overline{F}(x)\overline{G}(y),
\end{align}
where we used (\ref{lemma2res}) in the second step and Lemma \ref{lemma1} for $\gamma=0$. On the other hand, 
\begin{align}\label{Szuddd}
	&\quad P\left(X_1>ux,S_n>x, T_n>y\right)\nonumber\\
	&\geq  P\left(S_n>x, T_n>y\right)-P\left(S_n^{+}>x, T_n^{+}>y,X_1\leq ux \right)\nonumber\\
	&\geq P\left(S_n>x, T_n>y\right)-P\left(S_n^{+}>x, T_n^{+}>y\right)+P\left(S_n^{+}>x, T_n^{+}>y, X_1^{+}>x \right)\nonumber\\
	&=P\left(S_n>x, T_n>y\right)-P\left(S_n^{+}>x, T_n^{+}>y\right)+P\left(X_1^{+}>x, T_n^{+}>y\right),
\end{align}
where $S_n^{+}=\sum_{i=1}^nX_i^{+}$ and $T_n^{+}=\sum_{j=1}^nT_j^{+}$. It is straightforward to verify that if the pair $(X, Y)$ is distributed according to the FGM distribution, then the pair $(X^{+}, Y^{+})$ similarly follows the FGM distribution. Thus, it is obvious by (\ref{lemma1res}) that in the case of $\gamma=0$
\begin{align}\label{Szuddd1}
	P\left(S_n>x, T_n>y\right)\sim P\left(S_n^{+}>x, T_n^{+}>y\right)\sim n(n+\theta)\overline{F}(x)\overline{G}(y).
\end{align}
Since $(X_1^{+}, Y_1^{+})$ follows the FGM distribution, which satisfies the strongly asymptotic independence. Then, following Lemma 3.4 in Cheng et al. (2020), it can be derived that
\begin{align}\label{Szuddd2}
	P\left(X_1^{+}>x, T_n^{+}>y\right)\sim (n+\theta)\overline{F}(x)\overline{G}(y).
\end{align}
Substituting (\ref{Szuddd1}) and (\ref{Szuddd2}) into (\ref{Szuddd}) the reveals that
\begin{align*}
	P\left(X_1>ux,S_n>x, T_n>y\right)\gtrsim (n+\theta)\overline{F}(x)\overline{G}(y),
\end{align*}
which, together with (\ref{Szudd}), implies that (\ref{lemma2res2}) holds uniformly for $\varepsilon\leq u\leq 1$. As for $u>1$, it follows from (\ref{Szuddd2}) and the preceding result with $u=1$ respectively that
\begin{align*}
	P\left(X_1>ux,S_n>x, T_n>y\right)\leq P\left(X_1^{+}>ux, T_n^{+}>y\right)\sim (n+\theta)\overline{F}(ux)\overline{G}(y)
\end{align*}
and 
\begin{align*}
	P\left(X_1>ux,S_n>x, T_n>y\right)\geq P\left(X_1>ux,S_n>ux, T_n>y\right)\sim (n+\theta)\overline{F}(ux)\overline{G}(y).
\end{align*}
This completes the proof of (\ref{lemma2res2}).

Finally, we prove (\ref{lemma2res3}). For $u\geq \varepsilon$, the probability can be bounded from above in the following way.
\begin{align}\label{fenL}
	P\left(X_1^{-}>ux,S_n>x, T_n>y\right)&\leq \sum_{i=2}^nP\left(X_1^{-}>ux,X_i>\frac{x}{n}, Y_i>\frac{y}{n}\right)\nonumber\\
	&\quad +\sum_{i=2}^n\sum_{j=2,j\neq i}^nP\left(X_1^{-}>ux,X_i>\frac{x}{n}, Y_j>\frac{y}{n}\right)\nonumber\\
	&\quad +\sum_{i=2}^nP\left(X_1^{-}>ux,X_i>\frac{x}{n},Y_1>\frac{y}{n}\right)\nonumber\\
	&\triangleq L_1(x,y)+L_2(x,y)+L_3(x,y).
\end{align}
Since $F, G\in\mathscr{R}_{-\alpha}$, and $F(-x)=O\left(\overline{F}(x)\right)$, we have
\begin{align}\label{L1}
	L_1(x,y)&\lesssim (1+\theta)(n-1)n^{2\alpha}F(-ux)\overline{F}(x)\overline{G}(y)\nonumber\\
	&=O(1)\overline{F}(x)\overline{F}(ux)\overline{G}(y)\nonumber\\
	&=o(1)\overline{F}((u\vee 1)x)\overline{G}(y)
\end{align}
holds uniformly for $u>\varepsilon$. It is obvious that
\begin{align}\label{L2}
	L_2(x,y)=o(1)\overline{F}((u\vee 1)x)\overline{G}(y).
\end{align}
uniformly for $u>\varepsilon$. By utilizing (\ref{joint}), $F,G\in\mathscr{R}_{-\alpha}$, and $F(-x)=O\left(\overline{F}(x)\right)$ we derive that
\begin{align}\label{L3}
	L_3(x,y)&=O(1)\overline{F}(x)\left(F(-ux)-P\left(X_1\leq -ux,Y_1\leq \frac{y}{n}\right)\right)\nonumber\\
	&=O(1)\overline{F}(x)\left(F(-ux)-F(-ux)G\left(\frac{y}{n}\right)\left(1+\theta\overline{F}(-ux)\overline{G}\left(\frac{y}{n}\right)\right)\right)\nonumber\\
	&=O(1)\overline{F}(x)F(-ux)\overline{G}\left(\frac{y}{n}\right)\left(1-\theta \overline{F}(-ux)G\left(\frac{y}{n}\right)\right)\nonumber\\
	&=o(1)\overline{F}((u\vee 1)x)\overline{G}(y).
\end{align}
Plugging (\ref{L1})-(\ref{L3}) into (\ref{fenL}) yields (\ref{lemma2res3}). This completes the proof. $\hfill\blacksquare$

Building on the above lemmas, we will prove the main results in the remainder of this paper.
\subsection{Proof of Theorem \ref{thm1}}
\subsubsection*{$\bullet$ Proof of (\ref{thm1res})}
For any fixed $0< \varepsilon<\frac1{2}\wedge {\zeta}^{\beta}\wedge {(1-\zeta^{\beta})}$, we have
\begin{align}\label{thmdid}
 &\quad\rho_k(x,y;\beta, \zeta)\nonumber\\
	&=\left(\int_{[0,\varepsilon x^{\beta})\cup [(1-\varepsilon)x^{\beta}, x^{\beta})}+\int_{\varepsilon x^{\beta}}^{(1-\varepsilon) x^{\beta}}+\int_{x^{\beta}}^{\infty}\right)P\left(X_k>t^{\frac1{\beta}}\big\vert X_k>\zeta x, S_n>x, T_n>y\right)\mathrm{d}t\nonumber\\
	&\triangleq I_1(x,y)+I_2(x,y)+I_3(x,y).
\end{align}
It is easy to verify that
\begin{align}\label{thmI1}
	0\leq I_1(x,y)\leq 2\varepsilon x^{\beta}.
\end{align}
According to the proof of (\ref{lemma2res}), we notice that $$P\left(X_k>x(t^{\frac1{\beta}}\vee \zeta), S_n>x, T_n>y\right)\sim \frac1{n}P\left(S_n>x, T_n>y\right)$$ holds uniformly for $t\in[\varepsilon,1-\varepsilon)$. Then, by utilizing Lemma \ref{lemma2}, we have
\begin{align}\label{thmI2}
	I_2(x,y)&=x^{\beta}\int_{\varepsilon}^{1-\varepsilon}\frac{P\left(X_k>x(t^{\frac1{\beta}}\vee \zeta),S_n>x, T_n>y\right)}{P\left(X_k>\zeta x,S_n>x, T_n>y\right)}\mathrm{d}t\nonumber\\
	&\sim (1-2\varepsilon)x^{\beta}.
\end{align}
It follows from Lemma \ref{lemma2} and Remark \ref{re4.2} that
\begin{align}\label{thmI3}
	I_3(x,y)&=O(1)\int_{x^{\beta}}^{\infty}\frac{P\left(X_k>t^{\frac1{\beta}},T_n>y\right)}{\overline{F}(x)\overline{G}(y)}\mathrm{d}t\nonumber\\
	&=O(1)\int_{x^{\beta}}^{\infty}\frac{\overline{F}\left(t^{\frac1{\beta}}\right)}{\overline{F}(x)}\mathrm{d}t=O(1)x^{\beta}\int_{x}^{\infty}\frac{t^{\beta-1}\overline{F}\left(t\right)}{x^{\beta}\overline{F}(x)}\mathrm{d}t=o\left(x^{\beta}\right),
\end{align}
%\begin{align}\label{thmI3did}
%	I_3(x,y)&\lesssim n\int_{x^{\beta}}^{\infty}\frac{P\left(X_k>t^{\frac1{\beta}},T_n>y\right)}{\overline{F}(x)\overline{G}(y)}\mathrm{d}t\nonumber\\
%	&=n(1+\theta)\int_{x^{\beta}}^{\infty}\frac{P\left(X^{*}>t^{\frac1{\beta}}\right)P\left(T_{n-1}+Y^{*}>y\right)}{\overline{F}(x)\overline{G}(y)}\mathrm{d}t\nonumber\\
%	&\quad -n\theta  \int_{x^{\beta}}^{\infty}\frac{P\left(X_{\wedge}^{*}>t^{\frac1{\beta}}\right)P\left(T_{n-1}+Y^{*}>y\right)}{\overline{F}(x)\overline{G}(y)}\mathrm{d}t\nonumber\\
%	&\quad -n\theta \int_{x^{\beta}}^{\infty}\frac{P\left(X^{*}>t^{\frac1{\beta}}\right)P\left(T_{n-1}+Y_{\wedge}^{*}>y\right)}{\overline{F}(x)\overline{G}(y)}\mathrm{d}t\nonumber\\
%	&\quad +n\theta \int_{x^{\beta}}^{\infty}\frac{P\left(X_{\wedge}^{*}>t^{\frac1{\beta}}\right)P\left(T_{n-1}+Y_{\wedge}^{*}>y\right)}{\overline{F}(x)\overline{G}(y)}\mathrm{d}t\nonumber\\
%	&\triangleq n(1+\theta)I_{31}(x,y)-n\theta I_{32}(x,y)-n\theta I_{33}(x,y)+n\theta I_{34}(x,y).
%\end{align}
%The following relations can be verified directly by applying the change of variables and Lemma \ref{lemma-2}.
%\begin{align}\label{I31}
%	I_{31}(x,y)&\sim n\left(\widehat{G}(\gamma)\right)^{n-1}\frac{\beta\int_x^{\infty}t^{\beta-1}\overline{F}(t)\mathrm{d}t}{x^{\beta}\overline{F}(x)}x^{\beta}\nonumber\\
%	&=o\left(x^{\beta}\right),
%\end{align}
where the last step holds because of $F\in\mathscr{S}(\gamma)\subset \mathscr{R}_{-\infty}$, for $\gamma>0$, and Theorem A3.12(a) in {Embrechts et al. (1997)}. 
%Analogously, by exploiting some obvious tricks, one can get that 
%\begin{align}\label{I323334}
%	I_{32}(x,y)=o\left(x^{\beta}\right), \quad I_{33}(x,y)=o\left(x^{\beta}\right),\quad and \quad I_{34}(x,y)=o\left(x^{\beta}\right).
%\end{align}
%Plugging (\ref{I31}) and (\ref{I323334}) into (\ref{thmI3did}) yields that
%\begin{align}\label{thmI3}
%	0\leq I_3(x,y)=o\left(x^{\beta}\right).
%\end{align}
Then (\ref{thm1res}) is consequently proved by substituting (\ref{thmI1})-(\ref{thmI3}) into (\ref{thmdid}) and the arbitrariness of $\varepsilon$. 
\subsubsection*{$\bullet$ Proof of (\ref{thm1res-1})}
It can be expanded by the binomial theorem as follows.
\begin{align}\label{qingfen}
	\rho_k(x,y;\beta, 0)&=E\left(\sum_{i=0}^{\beta}(-1)^{\beta-i}\binom{\beta}{i}(X_k^{+})^{i}(X_k^{-})^{\beta-i}\Bigg\vert S_n>x, T_n>y\right)\nonumber\\
	&=\rho_k^{+}(x,y;\beta, 0)+(-1)^{\beta}\rho_k^{-}(x,y;\beta,0),
	%&=E\left((X_k^{+})^{\beta}\big\vert X_k>\zeta x, S_n>x, T_n>y\right)+(-1)^{\beta}E\left((X_k^{-})^{\beta}\big\vert X_k>\zeta x, S_n>x, T_n>y\right),
\end{align}
where $$\rho_k^{+}(x,y;\beta, 0)=E\left((X_k^{+})^{\beta}\big\vert S_n>x, T_n>y\right)$$ and $$\rho_k^{-}(x,y;\beta, 0)=E\left((X_k^{-})^{\beta}\big\vert S_n>x, T_n>y\right).$$
Through analogous manipulations, we conclude that
\begin{align}\label{qing+}
	\rho_k^{+}(x,y;\beta, 0)\sim \frac{x^{\beta}}{n}.
\end{align}
For some fixed $\varepsilon>0$, we derive that
\begin{align}\label{rho-}
	\rho_k^{-}(x,y;\beta, 0)&=\left(\int_{0}^{\varepsilon x^{\beta}}+\int_{\varepsilon x^{\beta}}^{\infty}\right)P\left(X_k^{-}>t^{\frac1{\beta}}\big\vert S_n>x, T_n>y\right)\mathrm{d}t\nonumber\\
	&\leq \varepsilon x^{\beta}+x^{\beta}\int_{\varepsilon}^{\infty}\frac{P\left(X_k<-t^{\frac1{\beta}}x, S_n>x, T_n>y\right)}{P\left(S_n>x, T_n>y\right)}\mathrm{d}t\nonumber\\
	&\leq \varepsilon x^{\beta}+x^{\beta}\int_{\varepsilon}^{\infty}\frac{P\left(\sum_{i=1,i\neq k}^nX_i>(1+t^{\frac1{\beta}})x, T_n>y\right)}{P\left(S_n>x, T_n>y\right)}\mathrm{d}t.
\end{align}
%Lemma \ref{lemma1}, (\ref{re4.22}) is applied in the first step and
Due to (\ref{re4.22}) and $F\in\mathscr{R}_{-\infty}$, we have for each fixed $t>0$ that
\begin{align*} 
	\frac{P\left(\sum_{i=1,i\neq k}^nX_i>\left(1+t^{\frac1{\beta}}\right)x, T_n>y\right)}{P\left(S_n>x, T_n>y\right)}=O(1) \frac{\overline{F}\left(\left(1+t^{\frac1{\beta}}\right)x\right)}{\overline{F}(x)}\to 0,\quad {\rm as}\quad x\to\infty.
\end{align*}
Furthermore, for $0<\delta<\gamma$ and $x$, $y$ large enough,
\begin{align*}
&\quad\frac{P\left(\sum_{i=1,i\neq k}^nX_i>\left(1+t^{\frac1{\beta}}\right)x, T_n>y\right)}{P\left(S_n>x, T_n>y\right)}\\
%	&=\frac{P\left(\sum_{i=1,i\neq k}^nX_i>(1+t^{\frac1{\beta}})x, T_n>y\right)}{P\left(\sum_{i=1,i\neq k}^nX_i>x, T_n>y\right)}\cdot \frac{P\left(\sum_{i=1,i\neq k}^nX_i>x, T_n>y\right)}{P\left(S_n>x, T_n>y\right)}\\
	&=O(1) \frac{\overline{F}\left(\left(1+t^{\frac1{\beta}}\right)x\right)}{\overline{F}(x)}=O(1)\frac{\overline{F}\left(x+\frac{t^{\frac1{\beta}}}{\gamma}\right)}{\overline{F}(x)}=O(1)\exp\left\{-\left(1-\frac{\delta}{\gamma}\right)t^{\frac{1}{\beta}}\right\},
\end{align*}
where (\ref{impor}) is adopted in the last step. It is easy to verify that $\exp\left\{-\left(1-\frac{\delta}{\gamma}\right)t^{\frac{1}{\beta}}\right\}$ is integrable with respect to $\mathrm{d}t$ over $(0,\infty)$. Then, by the dominated convergence theorem and the arbitrariness of $\varepsilon$, we arrive at
\begin{align}\label{qing-}
	\rho_k^{-}(x,y;\beta, 0)=o\left(\frac{x^{\beta}}{n}\right).
\end{align}
Substituting (\ref{qing+}) and (\ref{qing-}) into (\ref{qingfen}), we conclude (\ref{thm1res-1}). This closes the proof of Theorem \ref{thm1}.
$\hfill\blacksquare$
\subsection{Proof of Theorem \ref{thm2}}
\subsubsection*{$\bullet$ Proof of (\ref{thm2res})}
For any fixed $0<\zeta\leq 1$, we have 
\begin{align}\label{thm2did}
	\rho_k(x,y;\beta, \zeta)&=\left(\int_0^{(\zeta x)^{\beta}}+\int_{(\zeta x)^{\beta}}^{\infty}\right)P\left(X_k>t^{\frac1{\beta}}\big\vert X_k>\zeta x, S_n>x, T_n>y\right)\mathrm{d}t\nonumber\\
    &\triangleq T_1(x, y)+T_2(x,y).
\end{align}
It is easy to verify that 
\begin{align}\label{thm2T1}
	 T_1(x, y)=\zeta^{\beta}x^{\beta}.
\end{align}
By utilizing Lemma \ref{lemma2} and the uniform convergence for regular variation, we can get that
\begin{align}\label{thm2T2}
	T_2(x,y)&=x^{\beta}\int_{\zeta}^{\infty}\frac{\beta t^{\beta-1}P\left(X_k>tx, S_n>x, T_n>y\right)}{P\left(X_k>\zeta x, S_n>x, T_n>y\right)}\mathrm{d}t\nonumber\\
	&\sim x^{\beta}\int_{\zeta}^{\infty}\beta t^{\beta-1}(t\vee 1)^{-\alpha}\mathrm{d}t\nonumber\\
	&=\left(1-\zeta^{\beta}+\frac{\beta}{\alpha-\beta}\right)x^{\beta}.
\end{align}
Consequently, (\ref{thm2res}) is derived by substituting (\ref{thm2T1}) and (\ref{thm2T2}) into (\ref{thm2did}).
\subsubsection*{$\bullet$ Proof of (\ref{thm2res-1})}
Continuing with the notations used in the proof of (\ref{thm1res-1}), it can be verified that relation (\ref{qingfen}) still holds in the discussion of this case. Furthermore, we have
\begin{align}\label{zhong+}
	\rho_k^{+}(x,y;\beta, 0)\sim \frac{\alpha}{\alpha-\beta}\frac{x^{\beta}}{n}.
\end{align}
The proof of (\ref{zhong+}) is constructed in a manner analogous to the one used for the case of $0<\zeta\leq 1$. All that remains to prove Theorem \ref{thm2} is to prove (\ref{qing-}). According to (\ref{rho-})
\begin{align*}
	\rho_k^{-}(x,y;\beta, 0)&=\left(\int_{0}^{\varepsilon x^{\beta}}+\int_{\varepsilon x^{\beta}}^{\infty}\right)P\left(X_k^{-}>t^{\frac1{\beta}}\big\vert S_n>x, T_n>y\right)\mathrm{d}t\nonumber\\
	&\leq \varepsilon x^{\beta}+x^{\beta}\int_{\varepsilon}^{\infty}\frac{P\left(X_k^{-}>t^{\frac1{\beta}}x, S_n>x, T_n>y\right)}{P\left(S_n>x, T_n>y\right)}\mathrm{d}t\\
	&\sim \varepsilon x^{\beta}=o\left(\frac{x^{\beta}}{n}\right),\quad {\rm as}\quad \varepsilon\downarrow 0,
%	&=\left(\varepsilon +o(1)\int_{\varepsilon}^{\infty}\left(t\vee 1\right)^{-\frac{\alpha}{\beta}}\mathrm{d}t\right)x^{\beta}\\
%	&=o\left(\frac{x^{\beta}}{n}\right),\quad {\rm as}\quad \varepsilon\downarrow 0,
\end{align*}
where Lemmas \ref{lemma1} and \ref{lemma2} were utilized in the third setp. Thus, (\ref{qing-}) is proved, and we complete the proof of Theorem \ref{thm2}.
$\hfill\blacksquare$

\section*{Acknowledgment}
We are grateful to Prof. Yang Yang for his thoughtful suggestions and guidance, which greatly improved this work.

\end{document}